\documentclass{article}
\usepackage{amsmath}
\usepackage{makeidx}
\usepackage{amssymb}
\usepackage{amsfonts}

\newtheorem{theorem}{Theorem}
\newtheorem{acknowledgement}[theorem]{Acknowledgement}

\newtheorem{corollary}[theorem]{Corollary}

\newtheorem{definition}[theorem]{Definition}
\newtheorem{example}[theorem]{Example}

\newtheorem{lemma}[theorem]{Lemma}
\newtheorem{notation}[theorem]{Notation}
\newtheorem{problem}[theorem]{Problem}
\newtheorem{proposition}[theorem]{Proposition}
\newtheorem{remark}[theorem]{Remark}

\newenvironment{proof}[1][Proof]{\noindent\textbf{#1.} }{\ \rule{0.5em}{0.5em}}
\input{tcilatex}

\begin{document}

\title{Generalized Universal Covers of Uniform Spaces}
\author{Valera Berestovskii \\
Omsk Branch of the Sobolev Institute of Mathematics SD RAS \\
Pevtsova 13\\
Omsk\\
644099 Russia\\
berestov@iitam.omsk.net.ru \and Conrad Plaut \\
Department of Mathematics\\
University of Tennessee\\
Knoxville, TN 37919\\
USA\\
cplaut@math.utk.edu}
\maketitle
\date{}

\begin{abstract}
We develop a generalized covering space theory for a class of uniform spaces
called coverable spaces. Coverable spaces include all geodesic metric
spaces, connected and locally pathwise connected compact topological spaces,
in particular Peano continua, as well as more pathological spaces like the
topologist's sine curve. Each coverable space has a generalized universal
cover with universal and lifting properties. Associated with this
generalized universal cover is a functorial uniform space invariant called
the \textit{deck group}, which is related to the classical fundamental group
by a natural homomorphism. We obtain some specific results for
one-dimensional spaces.


\end{abstract}

\section{Introduction}

In this paper we construct a generalized universal cover for a very large
class of uniform spaces called \textit{coverable }spaces, which includes all
geodesic metric spaces (Corollary \ref{inner}), connected and locally
pathwise connected compact topological spaces (Corollary \ref{0628new}), and
in particular Peano continua (Corollary \ref{peano}). Coverable spaces also
include some more pathological spaces like the topologist's sine curve (see
below) and totally disconnected spaces (\cite{BPTG}). Associated with this
generalized universal cover is a functorial uniform space invariant called
the \textit{deck group}, which is related to the classical fundamental group
by a natural homomorphism.

Three types of uniform spaces are of greatest importance: topological
groups, compact topologial spaces and metric spaces. Topological groups were
considered in \cite{BPTG} and \cite{BPLCG}; the relationship to the present
paper may be found in Section \ref{topsec}, along with a correction to \cite%
{BPTG}. For compact topological spaces, which have a unique uniform
structure compatible with the topology, the deck group is a topological
invariant that coincides with the fundamental group in the setting of what
we will call Poincar\'{e} spaces (i.e., connected, locally pathwise
connected, semilocally simply connected spaces). The special case of metric
spaces is of particular interest. In addition to more classical examples
such as the Hawaiian earring and related spaces (\cite{G}, \cite{MM}, \cite%
{De}, \cite{CC2}, \cite{CC1}, \cite{B}), generalized universal covers and
fundamental groups have recently been studied in connection with
Gromov-Hausdorff limits of Riemannian manifolds (\cite{SW1} and \cite{SW2}).
Such limits are always geodesic spaces, and hence coverable. A recent
example of Menguy shows that limits of Riemannian manifolds with positive
Ricci curvature can have bad local topology (\cite{Mn})--precisely the sort
of metric spaces at which the present work is aimed. We will consider metric
spaces in more detail in an upcoming paper.

One of the main impediments to generalizing the classical construction of
the universal cover is the traditional definition of covering map, the most
important property of which is the ability to lift curves and homotopies.
However, this lifting property is traditionally gained at the expense of
requiring that a space and its cover be locally homeomorphic in a fairly
strong way, and consequently traditional universal (or simply connected)
covers exist \textit{only} for Poincar\'{e} spaces. Earlier work concerning
systems coverings of uniform spaces (cf. \cite{L}, \cite{D}) was also
limited by considering only traditional covers, and therefore a universal
object is impossible to obtain even for basic examples such as the countable
product of circles or the Hawaiian earring.

One may take a hint from topological groups for how to proceed. In this
category, a traditional covering map is a quotient homomophism with discrete
kernel. The action by the kernel is not only properly discontinuous, it is
uniformly so. Moreover, as we showed in \cite{BPTG}, one can exploit this
more uniform kind of action to define generalized covers as quotients with
central kernels that are (complete and) prodiscrete, i.e. are inverse limits
of discrete groups. These generalized covers have the lifting properties of
traditional covers, but are not in general local homeomorphisms. This allows
one to abandon assumptions concerning (semi-)local simple connectivity. The
kernel of the generalized universal cover is a kind of generalized
fundamental group.

The key to generalizing the results of \cite{BPTG} is to build a group
action right into the definition of generalized cover. Actually doing so is,
unfortunately, somewhat technical, and this was carried out by the second
author in \cite{P}. The basic idea is to consider a kind of
\textquotedblleft uniformly\textquotedblright\ properly discontinuous action
called a \textit{discrete action}\ on a uniform space. Inverse systems of
discrete actions define \textit{prodiscrete actions}\ that generalize the
notion of the action on a topological group via a prodiscrete, central
subgroup. The action by a subgroup also preserves the uniform structure (as
long as one matches the left or right action to the left or right
uniformity); this property is generalized by something called an
\textquotedblleft isomorphic\textquotedblright\ action on a uniform
structure, which also broadens the notion of an isometric action on a metric
space. A generalized \textit{cover} of uniform spaces is defined in \cite{P}
to be a quotient via a prodiscrete, isomorphic action. At the beginning of 
\cite{P} is a review of basic definitions and properties of uniform spaces
that are used for the construction; we will use the same notation in this
paper (see \cite{Bk} for a more in-depth discussion).

What follows is a sketch of our construction and main results. All of these
constructions involve a choice of basepoint(s), but basepoint choice has an
impact in this setting similar to that in the traditional setting for
pathwise connected spaces, and for simplicity we will save detailed
discussion of basepoints for the body of the paper. For now, all functions
are simply assumed to preserve some chosen basepoints (this is particularly
important for the uniqueness statements, which are only true up to choice of
basepoint).

For each uniform space $X$ there is an inverse system $(X_{E},\phi _{EF})$
indexed on the collection of entourages of $X$ partially ordered by reverse
inclusion, called the \textit{fundamental inverse system} of $X$. $X_{E}$
consists of equivalence classes of finite $E$-chains starting at a fixed
basepoint. An $E$-\textit{chain} between points $p,q\in X$ is an ordered set
of points $x_{0}:=p,...,x_{n}:=q$ such that for all $i$, $(x_{i},x_{i+1})\in
E$. An $E$\textit{-loop} is an $E$-chain that starts and ends at the same
point. Two $E$-chains from $p$ to $q$ are equivalent if one can be obtained
from the other through finitely many steps, each of which involves removing
or adding a point, always leaving the endpoints fixed and keeping an $E$%
-chain at each stage. The collection of all equivalence classes of $E$-loops
at the basepoint forms a group $\delta _{E}$ with respect to concatenation.
This group is finitely generated when $X$ is compact (Theorem \ref{compact}%
). When $X_{E}$ is provided with a natural \textquotedblleft
lifted\textquotedblright\ uniformity, $\delta _{E}$ acts discretely on $X_{E}
$ by concatenating an $E$-loop at the beginning of an $E$-chain. In a sense $%
\delta _{E}$ detects \textquotedblleft holes\textquotedblright\ in $X$ that
are, roughly speaking, \textquotedblleft larger than $E$\textquotedblright .
The restriction $\theta _{EF}$ of the bonding map $\phi _{EF}$ to $\delta
_{F}$ preserves concatenation and produces another inverse system $(\delta
_{E},\theta _{EF})$ of groups and homomorphisms; in fact these two systems
form an inverse system of actions as defined in \cite{P}. The inverse limit
of the group system consists of a group $\delta _{1}(X):=\underleftarrow{%
\lim }\delta _{E}$, called the \textit{deck group} of $X$, which acts
prodiscretely and isomorphically on $\widetilde{X}:=\underleftarrow{\lim }%
X_{E}$. Observe that $\delta _{1}(X)$ is actually a (prodiscrete)
topological group and the homomorphisms induced on the deck group are
continuous homomorphisms. However, we do not know of examples of spaces
having deck groups that are abstractly, but not continuously, isomorphic.

In general the projection $\phi :\widetilde{X}\rightarrow X$ may not be
surjective, and hence not a cover. In fact $\widetilde{X}$ may be only a
single point even when $X$ is not. The next definition, which is central to
the paper, deals with this issue.

\begin{definition}
\label{coverdef}Let $X$ be a uniform space. An entourage $E$ such that the
projection $\phi _{E}:\widetilde{X}\rightarrow X_{E}$ is surjective is
called a covering entourage. A uniform space $X$ is called coverable if
there is a uniformity base of covering entourages including $X\times X$. The
collection of all such entourages is called the covering base $\mathcal{C}%
(X) $.
\end{definition}

If $X$ is coverable then each $\phi _{E}$ is a covering map (in the sense of 
\cite{P}) and we refer to the projection $\phi :\widetilde{X}\rightarrow X$
as the \textit{universal cover of }$X$. The next theorem combines results of
Theorems \ref{73} and \ref{101}.

\begin{theorem}
(Induced Mapping) Let $X,Y$ be coverable spaces with universal covers $\phi :%
\widetilde{X}\rightarrow X$ and $\psi :\widetilde{Y}\rightarrow Y$. If $%
f:X\rightarrow Y$ is uniformly continuous then there is a unique uniformly
continuous function $\widetilde{f}:\widetilde{X}\rightarrow \widetilde{Y}$
such that $f\circ \phi =\psi \circ \widetilde{f}$. Moreover,

\begin{enumerate}
\item If $f$ is a cover then $\widetilde{f}$ is a uniform homeomorphism.

\item If $Z$ is a coverable space and $g:Y\rightarrow Z$ is uniformly
continuous then $\widetilde{g\circ f}=\widetilde{g}\circ \widetilde{f}$.
\end{enumerate}
\end{theorem}

The restriction of $\widetilde{f}$ to $\delta _{1}(X)$ in the above theorem
is a homomorphism $f_{\ast }:\delta _{1}(X)\rightarrow \delta _{1}(Y)$
(Theorem \ref{73}). Therefore the deck group is a functorial invariant of
uniform structures. In the case of compact spaces the deck group is a
topological invariant; if $X$ is a compact Poincar\'{e} space then the deck
group is naturally isomorphic to the fundamental group (Corollary \ref%
{0518equiv}). For non-compact spaces the deck group need not be a
topological invariant; for example the surface of revolution $S$ obtained by
rotating the graph of $e^{x}$ about the $x$-axis has trivial deck group even
though it is homeomorphic to a standard cylinder--the deck group of which is 
$\mathbb{Z}$. The problem is that a generator of $\pi _{1}(S)$ may be
represented by a path that extends down the cusp arbitrarily far, wraps
around the small cusp, and then travels back to the basepoint. This
generator will therefore not be detected by any of the groups $\delta
_{E_{\varepsilon }}$, where $E_{\varepsilon }$ is the metric entourage of
size $\varepsilon $. The deck group does indicate that the standard cylinder
and $S$ are not \textit{uniformly} homeomorphic and that $S$ is not
uniformly semilocally simply connected. The next theorem follows from
Corollary \ref{partialstep} and Theorem \ref{101}.

\begin{theorem}
(Universal Property) If $X$ and $Y$ are coverable and $f:X\rightarrow Y$ is
a cover then there is a unique cover $f_{B}:\widetilde{Y}\rightarrow X$ such
that $f\circ f_{B}=\phi $, where $\phi $ is the universal cover of $Y$.
\end{theorem}

With regard to the universal property it should be pointed out that, in
contrast to the situation for coverable topological groups, we do not know
whether the composition of covers between coverable spaces (or uniform
spaces in general) is a cover. (The situation for topological groups is
significantly simpler because the deck group is actually a central subgroup
of the universal cover and covers are simply quotient homomorphisms, making
their composition easier to understand.) Recall that the composition of
traditional covering maps between connected topological spaces need not be a
traditional cover (\cite{M}). We resolve this problem in the same way that
it is resolved for topological spaces (see \cite{S}): define a category
whose objects are covers $p:Y\rightarrow X$ between coverable uniform spaces
and whose morphisms are commutative diagrams 
\begin{equation*}
\begin{array}{lll}
X_{1} & \overset{f}{\longrightarrow } & \text{ \ \ \ \ \ \ }X_{2} \\ 
\text{ \ \ \ }\searrow ^{p_{1}} &  & ^{p_{2}}\swarrow  \\ 
& X & 
\end{array}%
\end{equation*}%
where $p_{1}$ and $p_{2}$ are covers and $f$ is uniformly continuous. It is
an immediate consequence of the preceding theorem that the universal cover $%
\phi :\widetilde{X}\rightarrow X$ is a universal object in this category.

\begin{definition}
\label{universaldef}A uniform space $X$ is called universal if there is a
base $\mathcal{U}$ for the uniformity such that for any $E\in \mathcal{U}$, $%
\phi _{XE}:X_{E}\rightarrow X$ is a uniform homeomorphism. The collection of
all such uniformities is called the universal base of $X$ (which always
contains $X\times X$).
\end{definition}

If $X$ is coverable then $\widetilde{X}$ is universal (Theorem \ref{xic})
and every universal space is coverable (Corollary \ref{uic}). Moreover, a
coverable space $Y$ is universal if and only if $\delta _{1}(X)$ is trivial,
or equivalently every $E$-loop is $E$-homotopic to the trivial loop for
entourages $E$ in a particular basis (Corollary \ref{uit}), a condition that
is reminiscent of simply connected. The next theorem follows from Theorem %
\ref{76}.

\begin{theorem}
(Lifting) Let $X$ be universal, $Y$ be coverable and $f:X\rightarrow Y$ be
uniformly continuous. Then there exists a unique uniformly continuous
function $f_{L}:X\rightarrow \widetilde{Y}$ such that $f=\phi \circ f_{L}$,
where $\phi $ is the universal cover of $Y$.
\end{theorem}

Any connected, uniformly locally pathwise connected (see Definition \ref%
{unidef}), simply connected uniform space is universal--in particular any
compact, connected, locally pathwise connected, simply connected topological
space is universal (Theorem \ref{0627simply} and Corollary \ref{compactu}).
Therefore we may apply the Lifting Theorem to paths and homotopies of paths.
This allows one to define a natural mapping $\lambda :\pi _{1}(X)\rightarrow
\delta _{1}(X)$, for any coverable space $X$, by lifting a loop that
represents an element of the fundamental group and taking the deck
transformation that takes the basepoint to the endpoint of the loop. Even
though the action is not discrete this mapping is well defined, and a
homomorphism.

There is a satisfying relationship between the two most basic algebraic
properties of the map $\lambda $ and topological properties of $\widetilde{X}
$ when $X$ is pathwise connected: (1) $\lambda $ is injective if and only if 
$\widetilde{X}$ is simply connected (Proposition \ref{inj}), and (2) $%
\lambda $ is surjective if and only if $\widetilde{X}$ is pathwise
connected. More precisely, for arbitrary coverable $X$, the image of $%
\lambda $ is the stabilizer in the deck group $\delta _{1}(X)$ of the
pathwise connected component of $\widetilde{X}$ that contains the basepoint
(Theorem \ref{nat}). Hence if $\widetilde{X}$ is both pathwise connected and
simply connected then the deck group and the fundamental group are
isomorphic. This is true as mentioned earlier when $X$ is a compact Poincar%
\'{e} space--and also when $X$ is a locally compact, pathwise connected
topological group (cf. \cite{BPLCG}). If $X$ is connected and uniformly
locally pathwise connected then the pathwise connected component of $%
\widetilde{X}$ is dense in $\widetilde{X}$ (Proposition \ref{dense}) and
therefore $\widetilde{X}$ is connected. In this case it is also true that $%
\lambda (\pi _{1}(X))$ is dense in $\delta _{1}(X)$. Note that, in the case
of the surface $S$ described above, $S$ is itself universal. $S$ is pathwise
connected but not simply connected and correspondingly $\lambda :\mathbb{%
Z=\pi }_{1}(X)\rightarrow \delta _{1}(X)=0$ is surjective but not injective.

As another example, consider the (closed) topologist's sine curve $T$,
illustrated in Figure \ref{tsc} with its universal cover. The deck group of $%
T$ is $\mathbb{Z}$; the action shifts the universal cover in a way similar
to the action of $\mathbb{Z}$ on $\mathbb{R}$. Note that arbitrarily fine
chains may wrap around $T$, while no path does. Hence the deck group
indicates topology where the fundamental group, which is trivial in this
case, does not. Finally, $\widetilde{T}$ is simply connected but not
pathwise connected, and correspondingly $\lambda :0=\pi _{1}(T)\rightarrow
\delta _{1}(T)=\mathbb{Z}$ is injective but not surjective (and the image of 
$\lambda $ need not be dense in the deck group because $T$ is not locally
pathwise connected).

\FRAME{ftbpFU}{3.5566in}{2.5413in}{0pt}{\Qcb{The Topologist's Sine Curve}}{}{%
topsine.jpg}{\special{language "Scientific Word";type "GRAPHIC";display
"USEDEF";valid_file "F";width 3.5566in;height 2.5413in;depth
0pt;original-width 10in;original-height 5.4296in;cropleft "0";croptop
"1";cropright "1";cropbottom "0";filename '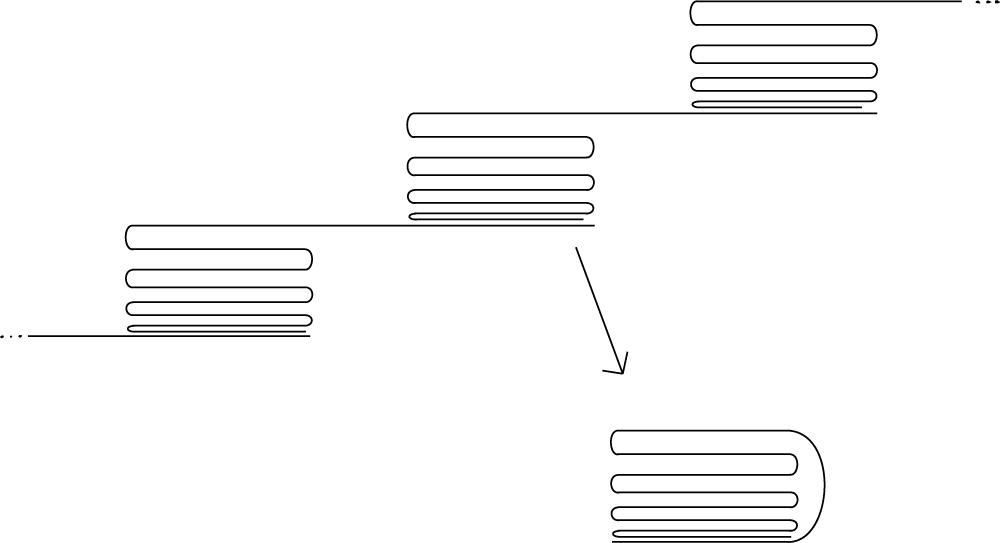';file-properties
"XNPEU";}}

We show that if $X$ is coverable with uniform dimension $\leq n$ in a sense
due to Isbell (\cite{I}) then $\widetilde{X}$ has the same property (Theorem %
\ref{dimthm}), and we conjecture that the the uniform dimensions of $X$ and $%
\widetilde{X}$ are the same. When $X$ has dimension $1$, we show that $%
\widetilde{X}$ is simply connected and contains no topological circles. As a
consequence the function $\lambda $ is injective. 

We prove some results concerning dimension and coverable spaces, although it
seems much more may be said. Further progress may require additional work
extending theorems concerning covering dimension to theorems about uniform
dimension. However, this partial result is particularly useful when $n=1$
because then $\widetilde{X}$ is forced to be at most one dimensional--and if 
$X$ contains any non-trivial curve then the uniform dimension of $\widetilde{%
X}$ must be exactly $1$.

One-dimensional metric spaces are of interest since they include, for
example, planar fractals and metrized Cayley graphs, as well as many
familiar pathalogical examples. It is known that the fundamental group of a
compact, connected, $1$-dimensional metrizable space embeds in an inverse
limit of free groups. Given this fact, it seems very likely that for
appropriate choice of $E$, the groups $\delta _{E}(X)$ are finitely
generated free groups when $X$ is compact, connected, and $1$-dimensional.
We have verified this in some special cases (see Section \ref{pssec}). It
would be very interesting to know whether for suitable choice of $E$, $%
\delta _{E}(X)$ is free when $X$ is coverable (or just connected?) of
uniform dimension $1$.

In this paper we show that, if $X$ is coverable of uniform dimension $\leq 1$
then $\widetilde{X}$ has no topological circles (Proposition \ref{circle}).
It follows that $\widetilde{X}$ is simply connected and hence (assuming $X$
is pathwise connected) $\lambda :\pi _{1}(X)\rightarrow \delta _{1}(X)$ is
an embedding. In general, $\widetilde{X}$ may not be pathwise connected and
hence $\lambda $ may not be surjective, but at least in the case when $X$ is
uniformly locally pathwise connected we know that the (isomorphic) image of
the fundamental group is dense in the deck group. In other words, $X$ is a
quotient via a free action of its fundamental group on a connected, simply
connected, uniformly one dimensional space.

We conclude with a general result about metric spaces. Given a metric space $%
X$, $\widetilde{X}$ is metrizable (Proposition \ref{metrizable}), and if $X$
is coverable then by definition of cover the action of $\delta _{1}(X)$ on $%
\widetilde{X}$ is isomorphic in the sense that there is a uniformity base
for $\widetilde{X}$ that is invariant with respect to $\delta _{1}(X)$. It
is natural to ask whether there is a metric on $\widetilde{X}$ with respect
to which the action of $\delta _{1}(X)$ is isometric and such that the
metric on $\widetilde{X}$ is the quotient metric. (Since $\delta _{1}(X)$ is
complete and acts prodiscretely, the orbits of $\delta _{1}(X)$ are closed
and hence if $\widetilde{X}$ has an invariant metric there is a well defined
quotient metric on $X$, namely the distance between the corresponding orbits
in $\widetilde{X}$--see \cite{Po}, 4.4 for details on quotient metrics.) We
do not know the answer to this question in general. However, if one examines
the proof in \cite{Bk} that a uniform space with a countable base possesses
a compatible pseudometric, it is clear that the explicitly constructed
pseudometric is invariant with respect to the action of a group $G$ if the
entourages used in the construction are invariant. In other words, if $X$ is
a coverable metric space then we may put a metric on $\widetilde{X}$ that is
invariant with respect to the action of $\delta _{1}(X)$. Although the
quotient metric on $X$ may not be the original, it is still uniformly
equivalent to the original. This proves:

\begin{theorem}
If $X$ is a coverable metric space then $X$ is uniformly homeomorphic to the
metric quotient $Y/G$ of a universal metric space $Y$ with respect to the
isometric prodiscrete action of a group $G$ isomorphic to $\delta _{1}(X)$.
\end{theorem}

\section{The Fundamental Inverse System}

Note that for any entourage $E$,%
\begin{equation}
E^{n}=\{(x_{0},x_{n})\in X\times X:\text{for some }x_{1},...,x_{n-1}\text{, }%
(x_{i},x_{i+1})\in E\text{ for }0\leq i<n\}\text{.}  \label{ton}
\end{equation}%
That is, $E^{n}$ is the set of all pairs of points joined by an $E$-chain of
length $n$.

\begin{definition}
A uniform space $X$ is called \textit{chain connected} (sometimes called
uniformly connected) if every uniformly continuous function from $X$ into a
discrete uniform space is constant.
\end{definition}

\begin{proposition}
\label{consume}The following are equivalent:

\begin{enumerate}
\item $X$ is chain connected

\item For any entourage $E$, $X\times X=\bigcup\limits_{n=1}^{\infty }E^{n}$.

\item For any entourage $E$, every pair of points in $X$ is joined by an $E$%
-chain.

\item For any entourage $E$ and $x\in X$, $X=\bigcup\limits_{n=1}^{\infty
}B(x,E^{n})$.
\end{enumerate}
\end{proposition}

\begin{proof}
The equivalence of the first two conditions is proved in \cite{J}, 9.34--the
statement of 9.34 is wrong but the proof is right! The last three are
equivalent by Formula (\ref{ton}).
\end{proof}

\begin{corollary}
If $X$ is chain connected and $f:X\rightarrow Y$ is a uniformly continuous
surjection then $Y$ is chain connected.
\end{corollary}

For a topological group $G$, chain connected is equivalent to $G$ being
locally generated (generated by every neighborhood of the identity), and $E$%
-chains are the same as what we called $U$-chains in \cite{BPTG}. If one
fixes a single point $p$ and entourage $E$ in an arbitrary uniform space $X$
then if $J_{p}\subset X$ is the set of all points that can be joined to $p$
by a $E$-chain, it is easy to check that $J_{p}$ is both open and closed.
Therefore every connected uniform space is chain connected. On the other
hand, the rational numbers are chain connected and totally disconnected (see 
\cite{BPTG} for related topics).

\begin{lemma}
\label{315system}Let $\{X_{\alpha },\phi _{\alpha \beta }\}_{\alpha \in
\Lambda }$ be an inverse system of sets such that each of the projections $%
\phi _{\alpha }:\underleftarrow{\lim }X_{\alpha }\rightarrow X_{\alpha }$ is
surjective. Then for any subsets $E,F$ of $X_{\alpha }\times X_{\alpha }$, $%
\phi _{\alpha }^{-1}(EF)=\phi _{\alpha }^{-1}(E)\phi _{\alpha }^{-1}(F)$. In
particular for any $n$, $\phi _{\alpha }^{-1}(E^{n})=\left( \phi _{\alpha
}^{-1}(E)\right) ^{n}$.
\end{lemma}

\begin{proof}
We have $((x_{\beta }),(y_{\beta }))\in \phi _{\alpha }^{-1}(EF)$ if and
only if there exists some $z_{\alpha }\in X_{\alpha }$ such that $(x_{\alpha
},z_{\alpha })\in E$ and $(z_{\alpha },y_{\alpha })\in F$. Since $\phi
_{\alpha }$ is surjective this is equivalent to: for some $(z_{\beta })\in X=%
\underleftarrow{\lim }X_{\beta }$, $((x_{\beta }),(z_{\beta }))\in \phi
_{\alpha }^{-1}(E)$ and $((y_{\beta }),(z_{\beta }))\in \phi _{\alpha
}^{-1}(F)$. But this is equivalent to $\left( (x_{\beta }\right) ,\left(
y_{\beta }\right) )\in \phi _{\alpha }^{-1}(E)\phi _{\alpha }^{-1}(F)$.
\end{proof}

\begin{lemma}
\label{41}Let $\{X_{\alpha },\phi _{\alpha \beta }\}_{\alpha \in \Lambda }$
be an inverse system of chain connected uniform spaces. If each of the
projections $\phi _{\alpha }$ is surjective then $X=\underleftarrow{\lim }%
X_{\alpha }$ is chain connected.
\end{lemma}

\begin{proof}
For any entourage $\phi _{\alpha }^{-1}(E)$ in $X$ and $\left( (x_{\beta
}),(y_{\beta })\right) \in X$ we have $\left( x_{\alpha },y_{\alpha }\right)
\in E^{n}$ for some $n$. Then $\left( (x_{\beta }),(y_{\beta })\right) \in
\phi _{\alpha }^{-1}((x_{\alpha },y_{\alpha }))\subset \phi _{\alpha
}^{-1}(E^{n})=\left( \phi _{\alpha }^{-1}(E)\right) ^{n}$.
\end{proof}

An $E$-extension of a $E$-chain $\{x_{0},...,x_{n}\}$ is an $E$-chain $%
\{x_{0},...,x_{i},x^{\prime },x_{i+1},...x_{n}\}$, where $0\leq i<n$. Two $E$%
-chains from $x_{0}$ to $x_{n}$ are said to be $E$-related if one is a $E$%
-extension of the other. An $E$-homotopy between $E$-chains $\gamma _{0}$
and $\gamma _{m}$ is a sequence $\{\gamma _{0},...,\gamma _{m}\}$ of $E$%
-chains such that $\gamma _{i}$ is $E$-related to $\gamma _{i-1}$ for all $%
1\leq i\leq m$. The number $m$ is called the \textit{length} of the
homotopy. We denote the $E$-homotopy class of an $E$-chain $\gamma $ by $%
[\gamma ]_{E}$. Now fix a basepoint $\ast \in X$. Let $X_{E}$ denote the set
of all $E$-homotopy classes of $E$-chains in $X$ starting at $\ast $ and
define a function (the \textquotedblleft endpoint map\textquotedblright ) $%
\phi _{XE}:X_{E}\rightarrow X$ by $\phi _{XE}([\ast =x_{0},...,x_{n}])=x_{n}$%
. The proof of the next lemma is immediate.

\begin{lemma}
\label{lemsur5}If $X$ is a uniform space and $E$ is an entourage then $\phi
_{XE}$ is surjective if and only if every pair in $X$ is joined by an $E$%
-chain; in particular $\phi _{XE}$ is surjective if $X$ is chain connected.
\end{lemma}

The next lemma will often be used without reference. The proof is
straightforward but tedious.

\begin{lemma}
\label{lemsur2}Let $\alpha :=\{\ast =a,x_{1},...,x_{n-1},b\}$ and $\beta
:=\{\ast =a,y_{1},...,y_{m-1},b\}$ be $E$-chains for some entourage $E$ in a
uniform space $X$. Then $\alpha $ and $\beta $ are $E$-homotopic if and only
if the $E$-loop $\alpha \ast \beta ^{-1}$ is $E$-homotopic to the trivial
chain $\{\ast \}$.
\end{lemma}

\begin{definition}
Let $X$ be a uniform space with entourage $E$. For any entourage $D\subset E$%
, define $D^{\ast }$ as follows: let $([\alpha ]_{E},[\beta ]_{E})\in
D^{\ast }$ if and only if 
\begin{equation*}
([\alpha ]_{E},[\beta ]_{E})=([\ast =x_{0},...,x_{n},y]_{E},[\ast
=x_{0},...,x_{n},z]_{E})\text{ with }(y,z)\in D\text{.}
\end{equation*}
\end{definition}

\begin{lemma}
\label{technical}Let $X$ be a uniform space with entourage $E$. For any
entourages $D,F\subset E$, $([\alpha ]_{E},[\beta ]_{E})\in D^{\ast }F^{\ast
}$ if and only if for some $[\gamma ]_{E}=[\ast =x_{0},...,x_{n}]_{E}$ we
have 
\begin{equation*}
\lbrack \alpha ]_{E}=[\ast =x_{0},...,x_{n},y]_{E}
\end{equation*}%
\begin{equation*}
\lbrack \beta ]_{E}=[\ast =x_{0},...,x_{n},z]_{E}
\end{equation*}%
for some $y,z$ such that $(y,x_{n})\in D$ and $(z,x_{n})\in F$. In
particular, $[\alpha ]_{E},[\beta ]_{E}\in B([\gamma ]_{E},\left( D^{\ast
}\right) ^{2})$ if and only if these conditions hold with $D=F$.
\end{lemma}

\begin{proof}
The reverse implication is obvious. Conversely, if $([\alpha ]_{E},[\beta
]_{E})\in D^{\ast }F^{\ast }$ then we may take $\alpha =\{\ast
=x_{0},...,x_{n-1},y\}$, $\gamma =\{\ast =x_{0},...,x_{n}\}$, $\gamma
^{\prime }=\{\ast =y_{0},...,y_{m},x_{n}\}$ and $\beta =\{\ast
=y_{0},...,y_{m},z\}$ where $\gamma ^{\prime }$ is $E$-homotopic to $\gamma $
and $(z,x_{n})\in F$ and $(y,x_{n})\in D$. Now $\beta $ is $E$-related to $%
\{\ast =y_{0},...,y_{m},x_{n},z\}$ and the latter is $E$-homotopic to $%
\{\ast =x_{0},...,x_{n},z\}$ using the $E$-homotopy from $\gamma $ to $%
\gamma ^{\prime }$. Finally, $\alpha $ is $E$-related to $\{\ast
=x_{0},...,x_{n},y\}$.
\end{proof}

\begin{proposition}
\label{xee}Let $X$ be a chain connected uniform space with entourage $E$,
and $\phi _{XE}:X_{E}\rightarrow X$ be the endpoint map. Then

\begin{enumerate}
\item $\phi _{XE}$ is injective when restricted to any $E^{\ast }$-ball.

\item For any entourage $D\subset E$ and $[\alpha ]_{E}\in X_{E}$, 
\begin{equation*}
\phi _{XE}(B([\alpha ]_{E},D^{\ast }))=B(\phi _{XE}([\alpha ]_{E}),D)
\end{equation*}%
and $\phi _{XE}(D^{\ast })=D$.

\item The collection of all $D^{\ast }$ such that $D\subset E$ is a base for
a uniformity of $X_{E}$.

\item $\phi _{XE}$ is a bi-uniformly continuous surjection with respect to
this uniformity. In particular the restriction of $\phi _{XE}$ to any $%
D^{\ast }$-ball is a uniform homeomorphism onto the corresponding $D$-ball .
\end{enumerate}
\end{proposition}

\begin{proof}
If $[\alpha ]_{E}$ and $[\beta ]_{E}$ are in $B([\gamma ]_{E},E^{\ast })$
and $\phi _{XE}([\alpha ]_{E})=\phi _{XE}([\beta ]_{E})$ then by definition
of $E^{\ast }$, $\alpha $ and $\beta $ are $E$-homotopic. This proves the
first part.

For Part (2) let $\alpha =\{\ast =x_{0},...,x_{n},x\}$; it is obvious from
the definition of $\phi _{XE}$ and $D^{\ast }$ that 
\begin{equation*}
\phi _{XE}(B([\alpha ]_{E},D^{\ast }))\subset B(\phi _{XE}([\alpha
]_{E}),D)=B(x,D)\text{.}
\end{equation*}%
Now suppose $(x,y)\in D$ and let $\beta =\{\ast =x_{0},...,x_{n},x,y\}$.
Since $\alpha $ is $E$-homotopic to $\{\ast =x_{0},...,x_{n},x,x\}$, we have
that $([\alpha ]_{E},[\beta ]_{E})\in D^{\ast }$ and 
\begin{equation*}
\phi _{XE}(([\alpha ]_{E},[\beta ]_{E}))=(x,y)\text{.}
\end{equation*}%
This implies that $\phi _{XE}(D^{\ast })=D$.

We will now check the conditions for a uniformity base. Clearly, for
entourages $E$ and $F$ we have $(E\cap F)^{\ast }\subset E^{\ast }\cap
F^{\ast }$ (and in fact they are equal but it is not necessary to prove
this). Since $D$ is symmetric, so is $D^{\ast }$. Next, let $F$ be an
entourage such that $F^{2}\subset D$. Suppose that $([\alpha ]_{E},[\beta
]_{E})\in (F^{\ast })^{2}$. Then for some $[\gamma ]_{E}$, $[\alpha
]_{E},[\beta ]_{E}\in B([\gamma ]_{E},F^{\ast })$. Applying Lemma \ref%
{technical} we may write $[\alpha ]_{E}=[\ast =x_{0},...,x_{n},y]_{E}$ and $%
[\beta ]_{E}=[\ast =x_{0},...,x_{n},z]_{E}$ with $(x_{n},z),(x_{n},y)\in F$
and hence $(y,z)\in F^{2}\subset D$. By definition $([\alpha ]_{E},[\beta
]_{E})\in D^{\ast }$. Recall from \cite{P} that $\phi _{XE}$ is bi-uniformly
continuous by definition if the image or inverse image of any entourage with
respect to $\phi _{XE}$ is again an entourage. Part (4) is now an immediate
consequence of Parts (2) and (3).
\end{proof}

We can now see precisely how the choice of basepoints affect things.

\begin{definition}
Let $\beta =\{x_{0},...,x_{n}\}$ and $\alpha =\{y_{0}=x_{n},...,y_{m}\}$ be $%
E$-chains for some entourage $E$ in a uniform space $X$. Define the
concatenation of $\alpha $ to $\beta $ by%
\begin{equation*}
\beta \ast \alpha :=\{x_{0},...,x_{n}=y_{0},...,y_{m}\}\text{. }
\end{equation*}
\end{definition}

\begin{remark}
\label{basepoint}Let $X$ be a chain connected uniform space, $p_{1},p_{2}\in
X$, and $E$ be an entourage. Let $X_{E}^{i}$ denote the space of $E$%
-homotopy classes of $E$-chains based at $p_{i}$. Let $\gamma $ be any $E$%
-chain from $p_{2}$ to $p_{1}$. There is a natural map from $X_{E}^{1}$ to $%
X_{E}^{2}$ defined by taking $[\alpha ]_{E}$ to $[\gamma \ast \alpha ]_{E}$,
where $\gamma \ast \alpha $ is the $E$-chain obtained by concatenating $%
\alpha $ to the end of $\gamma $. It is easy to check that this function is
a bijection with inverse $\eta \mapsto \gamma ^{-1}\ast \eta $. For any
entourage $F$ in $X$, there are corresponding entourages $F_{1}^{\ast }$ in $%
X_{E}^{1}$ and $F_{2}^{\ast }$ in $X_{E}^{2}$. Now $(\alpha ,\beta )\in
F_{1}^{\ast }$ if and only if $(\gamma \ast \alpha ,\gamma \ast \beta )\in
F_{2}^{\ast }$ and in particular, $X_{E}^{1}$ and $X_{E}^{2}$ are uniformly
homeomophic. Therefore, as in the case of traditional covering space theory,
the choice of basepoint\ $\ast $ plays only a minor and predictable role.
\end{remark}

\begin{notation}
Given a basepoint $\ast $ in $X$, we will always take $[\ast ]_{E}$ for the
basepoint in $X_{E}$. In general, for any function $f:X\rightarrow Y$ we
will always suppose that $f$ is \textquotedblleft pointed\textquotedblright\
in the sense that $f(\ast )=\ast $. For each chain $\gamma
=\{x_{0},...,x_{n}\}$ in $X$ we denote the chain $\{f(x_{0}),...,f(x_{n})\}$
by $f(\gamma )$, with similar notation for any finite sequence of chains. We
will always take the uniform structure on $X_{E}$ to be the one given by
Proposition \ref{xee}.
\end{notation}

If $E,F$ are entourages of $X,Y$, respectively, such that $f(E)\subset F$
and $\gamma $ is an $E$-chain then $f(\gamma )$ is an $F$-chain in $Y$. Note
that if $\eta $ is an $E$-homotopy between $E$-chains $\alpha $ and $\gamma $
then $f(\eta )$ is an $F$-homotopy between $f(\alpha )$ and $f(\gamma )$. In
particular the function in the following definition is well defined (and
pointed).

\begin{definition}
\label{induced}Given $f:X\rightarrow Y$ is a function between uniform spaces
such that for entourages $E,F$ of $X,Y$, respectively, $f(E)\subset F$, we
define $f_{EF}:X_{E}\rightarrow Y_{F}$ by $f([\gamma ]_{E})=[f(\gamma )]_{F}$%
.
\end{definition}

\begin{definition}
For any entourages $D\subset E$ in a uniform space $X$, define $\phi
_{ED}:X_{D}\rightarrow X_{E}$ by $\phi _{ED}([\alpha ]_{D})=[\alpha ]_{E}$.
\end{definition}

Note that by definition $\phi _{ED}=I_{DE}$, where $I:X\rightarrow X$ is the
identity, and 
\begin{equation*}
\phi _{XE}\circ \phi _{ED}=\phi _{XD}\text{.}
\end{equation*}

\begin{lemma}
For any entourages $D\subset E$ in a uniform space $X$, $\phi _{ED}$ is
uniformly continuous.
\end{lemma}

\begin{proof}
Let $F^{\ast }$ be an entourage in $X_{E}$ and $G\subset F\cap D$ be an
entourage in $X$. We will need to distinguish between $G^{\ast }\subset
X_{D}\times X_{D}$, which we will refer to as $G_{D}^{\ast }$ and $G^{\ast
}\subset X_{E}\times X_{E}$, which we will refer to as $G_{E}^{\ast }$. If $%
([\alpha ]_{D},[\beta ]_{D})\in G_{D}^{\ast }$ then by definition we can
take $\alpha =\{\ast =x_{0},...,x_{n},y\}$ and $\beta =\{\ast
=x_{0},...,x_{n},z\}$ with $(y,z)\in G$. Now 
\begin{equation*}
\phi _{ED}(([\alpha ]_{D},[\beta ]_{D}))=([\alpha ]_{E},[\beta ]_{E})
\end{equation*}%
and by definition $([\alpha ]_{E},[\beta ]_{E})\in G_{E}^{\ast }$. That is, $%
\phi _{ED}(G_{D}^{\ast })\subset G_{E}^{\ast }\subset F^{\ast }$.
\end{proof}

\begin{proposition}
\label{identify}Let $X$ be a uniform space and $D\subset E$ be entourages.
Then 
\begin{equation*}
\phi _{D^{\ast }D}:=(\phi _{XE})_{D^{\ast }D}:(X_{E})_{D^{\ast }}\rightarrow
X_{D}
\end{equation*}%
is a uniform homeomorphism such that the following diagram commutes:%
\begin{equation*}
\begin{array}{lll}
(X_{E})_{D^{\ast }} & \overset{\phi _{X_{E}D^{\ast }}}{\longrightarrow } & 
X_{E} \\ 
\downarrow ^{\phi _{D^{\ast }D}} & \nearrow _{\phi _{ED}} & \downarrow
^{\phi _{XE}} \\ 
X_{D} & \overset{\phi _{XD}}{\longrightarrow } & X%
\end{array}%
\end{equation*}%
Moreover, for any entourage $F\subset D$ in $X$, we have $\phi _{D^{\ast
}D}((F^{\ast })^{\ast })=F^{\ast }$.
\end{proposition}

\begin{proof}
For simplicity denote $(\phi _{XE})_{D^{\ast }D}$ by $\phi _{D^{\ast }D}$.
Let $C=\{[\ast ]_{E}=[\alpha _{0}]_{E},...,[\alpha _{n}]_{E}\}$ be a $%
D^{\ast }$-chain in $X_{E}$, where each $\alpha _{i}$ ends at a point $z_{i}$%
. Then 
\begin{equation}
\phi _{D^{\ast }D}([C]_{D^{\ast }})=[\phi _{XE}([\ast ]_{E}),...,\phi
_{XE}([\alpha _{n}]_{E})]_{D}=[z_{0}=\ast ,z_{1},...,z_{n}]_{D}\text{.}
\label{identifyfor}
\end{equation}

Since $([\ast ]_{E},[\alpha _{1}]_{E})\in D^{\ast }$, 
\begin{equation*}
\lbrack \alpha _{1}]_{E}=[z_{0}=\ast ,z_{1}]_{E}\text{.}
\end{equation*}%
Proceeding inductively, for all $i$ we have that 
\begin{equation}
\lbrack \alpha _{i}]_{E}=[z_{0}=\ast ,z_{1},...,z_{i}]_{E}
\label{0621sequence}
\end{equation}%
which implies that $\phi _{D^{\ast }D}$ is injective.

Given any $D$-chain $\{z_{0}=\ast ,z_{1},...,z_{n}\}$ we can let $\alpha
_{i}=\{\ast ,z_{1},..,z_{i}\}$ for all $0\leq i\leq n$. Since $\{z_{0}=\ast
,z_{1},...,z_{n}\}$ is a $D$-chain, $\{[\alpha _{0}]_{D},...,[\alpha
_{n}]_{D}\}$ is a $D^{\ast }$-chain and 
\begin{equation*}
\phi _{DD^{\ast }}([[\ast ]_{E}=[\alpha _{0}]_{E},...,[\alpha
_{n}]_{E}]_{D^{\ast }})=[z_{0}=\ast ,z_{1},...,z_{n}]_{D}
\end{equation*}%
which shows that $\phi _{D^{\ast }D}$ is surjective.

By definition, 
\begin{equation*}
\phi _{ED}\left( \phi _{D^{\ast }D}([C]_{D^{\ast }})\right) =[z_{0}=\ast
,z_{1},...,z_{n}]_{E}\text{.}
\end{equation*}%
That is, the upper triangle of the diagram commutes, and we already know
that the bottom triangle does. Finally, suppose that 
\begin{equation*}
([[\ast ]_{E}=[\alpha _{0}]_{E},...,[\alpha _{n}]_{E}]_{D^{\ast }},[[\ast
]_{E}=[\beta _{0}]_{E},...,[\beta _{m}]_{E}]_{D^{\ast }})\in (F^{\ast
})^{\ast }
\end{equation*}%
where, according to Formula \ref{0621sequence} we can suppose that for some $%
z_{1},...,z_{n},w_{1},...,w_{m}\in X$, $\alpha _{i}=\{\ast ,z_{1},..,z_{i}\}$
and $\beta _{i}=\{\ast ,w_{1},..,w_{i}\}$ for all $i$. By definition this
means that we can suppose that $m=n$ and $([\alpha _{n}]_{E},[\beta
_{n}]_{E})\in F^{\ast }$. This in turn is equivalent to $(z_{n},w_{n})\in F$%
, which is equivalent to 
\begin{equation*}
\phi _{DD^{\ast }}(([[\ast ]_{E}=[\alpha _{0}]_{E},...,[\alpha
_{n}]_{E}]_{D^{\ast }},[[\ast ]_{E}=[\beta _{0}]_{E},...,[\beta
_{m}]_{E}]_{D^{\ast }}))\in F^{\ast }\text{.}
\end{equation*}
\end{proof}

\begin{remark}
The preceding somewhat technical-looking proposition in fact has a very nice
interpretation. Essentially it identifies $(X_{E})_{D^{\ast }}$ with $X_{D}$
by taking a $D^{\ast }$-chain of $E$-chains to the $D$-chain of their
endpoints. In other words, $X_{E}$ and $X$ are \textquotedblleft locally the
same\textquotedblright ; $D^{\ast }$ and $E^{\ast }$ are really just copies
of $D$ and $E$ inside $X_{E}\times X_{E}$. Dealing with $X_{D}$ rather than $%
(X_{E})_{D^{\ast }}$ means we are dealing with chains rather than chains of
chains. At the same time, Proposition \ref{identify} identifies the mapping $%
\phi _{X_{E}D^{\ast }}:(X_{E})_{D^{\ast }}\rightarrow X_{E}$ with the more
easily understood mapping $\phi _{ED}:X_{D}\rightarrow X_{E}$. With such
identifications we can express the proposition as $(X_{E})_{D^{\ast }}=X_{D}$%
. This proposition is very useful because it will allow us to immediately
apply results proved for $\phi _{XE}$ to the more general functions $\phi
_{EF}$. The next lemma illustrates this.
\end{remark}

\begin{lemma}
\label{chainhom}Suppose that $X$ is a uniform space and $E$ is an entourage
in $X$ such that $X_{E}$ is chain connected. Then for any entourage $F$ and $%
E$-chain $\gamma $ in $X$, $\gamma $ is $E$-homotopic to an $F$-chain.
\end{lemma}

\begin{proof}
Without loss of generality we can suppose that $F\subset E$. Since $X_{E}$
is chain connected, $\phi _{X_{E}F^{\ast }}=\phi _{EF}$ is surjective by
Lemma \ref{lemsur5}. But this means that if $\gamma $ is an $E$-chain in $X$
there is some $F$-chain $\alpha $ in $X$ such that $[\alpha ]_{E}=\phi
_{EF}([\alpha ]_{F})=[\gamma ]_{E}$.
\end{proof}

Given entourages $D\subset E\subset F$ in $X$, we have functions $\phi
_{ED}:X_{D}\rightarrow X_{E}$ and $\phi _{FE}:X_{E}\rightarrow X_{F}$ with 
\begin{equation*}
\phi _{FE}(\phi _{ED}([x_{0}=\ast ,...,x_{n}]_{D}))=[x_{0}=\ast
,...,x_{n}]_{F}=\phi _{FD}([x_{0}=\ast ,...,x_{n}]_{D})\text{.}
\end{equation*}%
In other words, $\phi _{FE}\circ \phi _{ED}=\phi _{FD}$ and $\{X_{E},\phi
_{ED}\}$ forms an inverse system of uniformly continuous functions having as
its indexing set the set of all entourages of $X$ partially ordered by
reverse inclusion.

\begin{definition}
We will call the inverse system $\{X_{E},\phi _{ED}\}$ the fundamental
inverse system of $X$ and let $\widetilde{X}$ denote the inverse limit of
this inverse system with the inverse limit uniformity. We will let $\phi
_{E}:\widetilde{X}\rightarrow X_{E}$ be the projection; we simply denote $%
\phi _{X}$ by $\phi $. We will always choose for our basepoint $\ast $ in $%
\widetilde{X}$ the element having as each of its coordinates the basepoint $%
\ast =[\ast ]_{E}$ in $X_{E}$.
\end{definition}

Since a uniform space is metrizable if and only if it has a countable base
for its uniformity (\cite{Bk}), if $X$ is metrizable then each $X_{E}$ is
metrizable. We may index the fundamental system with this countable base and
conclude:

\begin{proposition}
\label{metrizable}If $X$ is metrizable then $\widetilde{X}$ is metrizable.
\end{proposition}

\begin{theorem}
\label{60}Let $X$ and $Y$ be uniform spaces, $f:X\rightarrow Y$ be uniformly
continuous, and $E,F$ be entourages in $X,Y$, respectively, such that $%
f(E)\subset F$. Then $f_{EF}$ is uniformly continuous and satisfies 
\begin{equation}
f\circ \phi _{XE}=\phi _{YF}\circ f_{EF}\text{.}  \label{commute1}
\end{equation}%
Moreover, if $X_{E}$ is chain connected then $f_{EF}$ is the unique
uniformly continuous function satisfying (\ref{commute1}) such that $%
f_{EF}([\ast ]_{E})=[f(\ast )]_{F}$.
\end{theorem}

\begin{proof}
To prove $f_{EF}$ is uniformly continuous we need only consider entourages
of the form $D^{\ast }$ in $Y_{F}$, where $D\subset F$ is an entourage in $Y$%
. Since $f$ is uniformly continuous there exists some entourage $G\subset E$
such that $f(G)\subset D$. Let $\alpha :=\{\ast =x_{0},...,x_{n},y\}$ and $%
\beta :=\{\ast =x_{0},...,x_{n},z\}$ be such that 
\begin{equation*}
([\alpha ]_{E},[\beta ]_{E})\in G^{\ast }\text{,}
\end{equation*}%
which means by definition $(y,z)\in G$. Then 
\begin{equation*}
(f(y),f(z))\in D
\end{equation*}%
and therefore 
\begin{equation*}
f_{EF}([\alpha ]_{E},[\beta ]_{F})=([\ast
=f(x_{0}),...,f(x_{n}),f(y)]_{F},[\ast =f(x_{0}),...,f(x_{n}),f(z)]_{F})\in
D^{\ast }\text{.}
\end{equation*}%
This shows that $f_{EF}$ is uniformly continuous. That $f\circ \phi
_{XE}=\phi _{YF}\circ f_{EF}$ is an immediate consequence of the definition.

To prove the last statement suppose $X_{E}$ is chain connected and $%
f^{\prime }:X_{E}\rightarrow Y_{F}$ is a uniformly continuous function such
that $f\circ \phi _{XE}=\phi _{YF}\circ f^{\prime }$ and $f^{\prime }([\ast
]_{E})=[f(\ast )]_{F}$. Let $G$ be an entourage in $X$ such that $f^{\prime
}(G^{\ast })\cup f_{EF}(G^{\ast })\subset F^{\ast }$. By way of Lemma \ref%
{chainhom} it is sufficient to show that if $\beta :=\{\ast
=x_{0},...,x_{n}\}$ is a $G^{\ast }$-chain then $f^{\prime }([\beta
]_{E})=f_{EF}([\beta ]_{F})$. We will prove it by induction on $n$. The case 
$n=0$ is given; suppose the statement is true for $n\geq 0$ and consider $%
\beta :=\{\ast =x_{0},...,x_{n+1}\}$ with $\alpha :=\{\ast =x_{0},...,x_{n}\}
$. Suppose that $f^{\prime }([\beta ]_{E}):=[f(\ast )=z_{0},...,z_{m}]_{F}$;
then 
\begin{equation*}
f(x_{n+1})=f\circ \phi _{XE}([\ast =x_{0},...,x_{n+1}]_{E})=\phi _{YF}\circ
f^{\prime }([\ast =x_{0},...,x_{n+1}]_{E})=z_{m}\text{.}
\end{equation*}%
By the inductive hypothesis, $f^{\prime }([\alpha ]_{E})=f_{EF}([\alpha
]_{E})$ and by definition of $G^{\ast }$ we have 
\begin{equation*}
(f_{EF}([\alpha ]_{E}),f^{\prime }([\beta ]_{E}))=(f^{\prime }([\alpha
]_{E}),f^{\prime }([\beta ]_{E}))\in F^{\ast }
\end{equation*}%
and 
\begin{equation*}
(f_{EF}([\alpha ]_{E}),f_{EF}([\beta ]_{E}))\in F^{\ast }\text{.}
\end{equation*}%
In other words, both $f_{EF}([\beta ]_{E})$ and $f^{\prime }([\beta ]_{E}))$
lie in $B(f_{EF}([\alpha ]_{E},F^{\ast })$ and by Proposition \ref{xee} $%
\phi _{YF}$ is injective on this ball. The fact that $f_{EF}([\beta
]_{E})=f^{\prime }([\beta ]_{E})$ now follows from 
\begin{equation*}
\phi _{YF}(f_{EF}([\beta ]_{E})=f(\phi _{XE}([\beta ]_{E}))=\phi
_{YF}(f^{\prime }([\beta ]_{E}))\text{.}
\end{equation*}
\end{proof}

\begin{corollary}
If $X$ is a uniform space and $F\subset E$ are entourages in $X$ such that $%
X_{F}$ is chain connected then $\phi _{EF}$ is the unique uniformly
continuous function such that $\phi _{EF}([\ast ]_{F})=[\ast ]_{E}$ and $%
\phi _{XF}=\phi _{XE}\circ \phi _{EF}$.
\end{corollary}

The proof of the next lemma is virtually identical to the proof of Lemma 65
in \cite{BPTG}; one need only replace statements like $xy^{-1}\in U$ with $%
(x,y)\in E$.

\begin{lemma}[Chain Lifting]
Let $X,Y$ be uniform spaces, $f:X\rightarrow Y$ be a uniformly continuous
surjection, $F$ be an entourage in $Y$ and $E:=f^{-1}(F)$. Let $c$ be an $E$%
-chain in $X$ and $\eta $ be an $F$-homotopy from the $F$-chain $d:=f(c)$ to
another $F$-chain $d^{\prime }$. Then $\eta $ lifts to an $E$-homotopy from $%
c$ to an $E$-chain $c^{\prime }$. That is, there exist an $E$-chain $%
c^{\prime }$ and an $E$-homotopy $\kappa $ between $c$ and $c^{\prime }$
such that $f(\kappa )=\eta $.
\end{lemma}

\begin{proposition}
\label{66}Let $X,Y$ be uniform spaces, $f:X\rightarrow Y$ be a uniformly
continuous surjection, $F$ be an entourage in $Y$ and $E:=f^{-1}(F)$. If $%
\phi _{XE}:X_{E}\rightarrow X$ is surjective and there exists a uniformly
continuous function $\psi :X\rightarrow Y_{F}$ such that the following
diagram commutes 
\begin{equation*}
\begin{array}{lll}
X_{E} & \overset{\phi _{XE}}{\longrightarrow } & X \\ 
\downarrow ^{f_{EF}} & \swarrow _{\psi } & \downarrow ^{f} \\ 
Y_{F} & \overset{\phi _{YF}}{\longrightarrow } & Y%
\end{array}%
\end{equation*}%
then $\phi _{XE}$ is a uniform homeomorphism.
\end{proposition}

\begin{proof}
We need only show that $\phi _{XE}$ is injective. Equivalently we need only
show that if $c:=\{\ast =x_{0},...,x_{n}=\ast \}$ is an $E$-loop in $X$ then 
$c$ is $E$-homotopic to the trivial loop $\{\ast \}$. Let $d:=f(c)$, which
is an $F$-loop in $Y$ and $f_{EF}([c]_{E}):=[h]_{F}$; by the commutativity
of the diagram $\phi _{YF}([h]_{F})=f(\phi _{XE}([c]_{E})=f(\ast )$. In
particular, $h$ is a loop in $Y$. Moreover, $[h]_{F}=\psi (\phi
_{XE}([c]_{E}))=\psi (\ast )$. On the other hand, 
\begin{equation*}
\lbrack f(\ast )]_{F}=f_{EF}([\ast ]_{E})=\psi (\phi _{XE}([\ast
]_{E}))=\psi (\ast )\text{. }
\end{equation*}%
That is, $[h]_{F}=[f(\ast )]_{F}$ and therefore $h$ is $F$-homotopic to the
trivial loop $f(\ast )$. Now by definition, 
\begin{equation*}
\lbrack h]_{F}=f_{EF}([c]_{E})=[f(x_{0}),...,f(x_{n})]_{F}=[f(c)]_{F}=[d]_{F}
\end{equation*}%
and therefore $d$ is also $F$-homotopic to $[f(\ast )]_{F}$. The Chain
Lifting Lemma now finishes the proof.
\end{proof}

\begin{corollary}
\label{67}Let $X$ and $Y$ be uniform spaces with $X$ chain connected, $%
f:X\rightarrow Y$ be a surjective uniformly continuous map, $F$ be an
entourage in $Y$ and $E:=f^{-1}(F)$. If $\phi _{YF}:Y_{F}\rightarrow Y$ is
bijective then $\phi _{XE}:X_{E}\rightarrow X$ is a uniform homeomorphism.
\end{corollary}

\begin{proof}
Note that since $X$ is chain connected, $\phi _{XE}$ is surjective. Let $%
\psi :=\phi _{YF}^{-1}\circ f$, which is uniformly continuous since $\phi
_{YF}$ is bi-uniformly continuous. Certainly $f=\phi _{YF}\circ \psi $ and 
\begin{equation*}
\phi _{YF}\circ (\psi \circ \phi _{XE})=f\circ \phi _{XE}=\phi _{YF}\circ
f_{EF}\text{.}
\end{equation*}%
Since $\phi _{YF}$ is bijective we conclude that $\psi \circ \phi
_{XE}=f_{EF}$ and the conditions of Proposition \ref{66} are satisfied,
completing the proof.
\end{proof}

\begin{proposition}
\label{coverchain}If $X$ is a coverable uniform space and $E$ is a covering
entourage in $X$ then $X_{E}$ is coverable and hence chain connected.
\end{proposition}

\begin{proof}
If $E$ is a covering entourage then we may index the fundamental inverse
system of $X_{E}$ using the set of all entourages $F^{\ast }$ where $%
F\subset E$ is a covering entourage in $X$. Then $\phi _{F}:\widetilde{X}%
\rightarrow X_{F}$ is surjective for all such $F$. But $(X_{E})_{F^{\ast }}$
is naturally identified with $X_{F}$ by Proposition \ref{identify} and
therefore $\phi _{F^{\ast }}:\widetilde{X_{E}}\rightarrow (X_{E})_{F^{\ast }}
$ is also surjective. Therefore the collection of all such $F^{\ast }$ is a
covering basis for $X_{E}$.
\end{proof}

It is straightforward but tedious to check that if $[\alpha ]_{E}=[\alpha
^{\prime }]_{E}$ and $[\beta ]_{E}=[\beta ^{\prime }]_{E}$ then $[\alpha
\ast \beta ]_{E}=[\alpha ^{\prime }\ast \beta ^{\prime }]_{E}$.

\begin{definition}
For an entourage $E$ in a uniform space $X$, we define $\delta _{E}(X)$ to
be the group of all $E$-homotopy classes of $E$-loops at $\ast $ with the
group operation induced by concatenation. We will call this the $E$-deck
group of $X$.
\end{definition}

That is, given $E$-loops $\alpha $ and $\beta $ based at $\ast $, we let 
\begin{equation*}
\lbrack \alpha ]_{E}\ast \lbrack \beta ]_{E}:=[\alpha \ast \beta ]_{E}\text{.%
}
\end{equation*}%
Note that the identity chain is $[\ast ]_{E}$ and if $\alpha =\{\ast
=x_{0},...,x_{n}=\ast \}$ then $[\alpha ]_{E}^{-1}=[\alpha ^{-1}]_{E}$ where 
$\alpha ^{-1}:=\{\ast =x_{n},...,x_{0}=\ast \}$. It is easy to check that $%
\delta _{E}(X)$ is in fact a group.

\begin{theorem}
\label{compact}If $X$ is a compact uniform space and $E$ is an entourage in $%
X$ such that $X_{E}$ is chain connected then $\delta _{E}(X)$ is finitely
generated.
\end{theorem}

\begin{proof}
Let $F$ be an entourage in $X$ such that $F^{3}\subset E$. Since $X$ is
compact there exists some finite $F$-dense set, i.e., a set $%
A:=\{x_{1},...,x_{n}\}$ such that for every $x\in X$ there exists some $%
x_{i}\in A$ such that $(x,x_{i})\in F$. We will first show that if $\alpha
:=\{y_{1},...,y_{m}\}$ is any $E$-chain with $y_{1},y_{m}\in A$ then $\alpha 
$ is $E$-homotopic to an $E$-chain $\{y_{1}=z_{1},z_{2},...,z_{m}=y_{m}\}$
such that $z_{i}\in A$ for all $i$. By Lemma \ref{chainhom} we may assume
that $\alpha $ is in fact an $F$-chain. Now for each $y_{i}$, $1<i<m$, there
is some $z_{i}$ such that $(y_{i},z_{i})\in F$. We may now proceed
iteratively, removing each $y_{i}$ and then replacing it by $z_{i}$. For
example, since $(y_{1},y_{3})\in F^{2}\subset E$, we may remove $y_{2}$ and
still have an $E$-chain. Since $(y_{1},y_{2})\in F$ and $(y_{2},z_{2})\in F$%
, $(y_{1},z_{2})\in F^{2}\subset E$. Similarly $(z_{2},y_{3})\in F^{2}$, so
we may add $z_{2}$. Likewise, $(z_{2},y_{4})\in F^{3}\subset E$, so we may
remove $y_{3}$, and then $(z_{3},z_{2})\in F^{3}$ and $(z_{3},y_{4})\in
F^{2} $, so we may add $z_{3}$. After finitely many steps we have the
desired $E$-chain.

We may suppose that $\ast \in A$. We now claim that $\delta _{E}(X)$ is
generated by elements of the form 
\begin{equation*}
\lbrack \ast =y_{1},...,y_{k},y_{k+1},...,y_{m}=y_{k},y_{k-1},...,y_{1}=\ast
]_{E}
\end{equation*}%
where $y_{i}\in A$ for all $i$, and for $1\leq i<j\leq m-1$ we have $%
y_{i}\neq y_{j}$. We will call elements of this form minimal elements. That
is, a minimal element is represented by a loop that consists of a chain made
of distinct points of $A$, followed by a loop (which may be empty) of
additional distinct points of $A$, followed by the initial chain in reverse
order. If we prove this claim then the proof of the theorem is finished
because $A$, and hence the set of minimal elements, is finite.

Let $\gamma =\{\ast =z_{1},...,z_{r}=\ast \}$ be an arbitrary $E$-loop. We
will show by induction on $r$ that $[\gamma ]_{E}$ is a product of minimal
elements. For $r=1$ the proof is trivial; suppose it is true for $r-1\geq 1$%
. If all of the points $z_{i}$ $1\leq i<r$ are distinct then $[\gamma ]_{E}$
is already a minimal element. Otherwise, let $j<r$ be the smallest index
such that for some $i<j$, $z_{i}=z_{j}$. Let 
\begin{equation*}
\beta :=\{z_{1},...,z_{i},z_{j+1},...,z_{r}\}\text{.}
\end{equation*}%
We may apply the inductive hypothesis to conclude that $\left[ \beta \right]
_{E}$ is the product of minimal elements. On the other hand, let 
\begin{equation*}
\alpha :=\{z_{1},...,z_{i},z_{i+1},...,z_{j}=z_{i},z_{i-1},...,z_{1}\}\text{.%
}
\end{equation*}%
By construction, $\alpha $ is minimal and since $\left[ \gamma \right] _{E}=%
\left[ \alpha \ast \beta \right] _{E}$, the proof is finished.
\end{proof}

The following mapping is certainly well-defined.

\begin{definition}
Let $X$ be a uniform space and $E$ be an entourage. For each $[\lambda
]_{E}\in \delta _{E}(X)$, define a mapping $\overline{\lambda }%
:X_{E}\rightarrow X_{E}$ by $\overline{\lambda }([\alpha ]_{E})=[\lambda
\ast \alpha ]_{E}$.
\end{definition}

Before proceeding we will recall some notation and definitions from \cite{P}%
. Let $X$ be a uniform space. We denote by $H_{X}$ the topological group of
uniform homeomorphisms of $X$ with composition as the operation. Suppose $G$
is a subgroup of $H_{X}$. An entourage $E$ is called $G$\textit{-invariant}
if $f(E)=E$ for every $f\in G$. If $X$ has a uniformity base consisting of $G
$-invariant entourages we say that $G$ acts isomorphically. The action of $G$
is said to be \textit{discrete} provided there exists some entourage $E$
such that if $(g(x),x)\in E$ for some $x\in X$ and $g\in G$ then $g$ must be
the identity. If $G$ acts discretely and isomorphically on $X$ then the
natural mapping $\pi :X\rightarrow X/G$ is called a discrete cover.

\begin{theorem}
\label{discrete}Let $X$ be a uniform space and $E$ be an entourage.

\begin{enumerate}
\item For any $\lambda _{1},\lambda _{2}$ we have $\overline{\lambda _{1}}%
\circ \overline{\lambda _{2}}=\overline{\lambda _{1}\ast \lambda _{2}}$ and $%
\delta _{E}(X)$ is naturally isomorphic to a subgroup of $H_{X_{E}}$.

\item $\delta _{E}(X)$ acts discretely and isomorphically on $X_{E}$.

\item If $\phi _{XE}$ is surjective (in particular if $X$ is chain
connected) then $\phi _{XE}:X_{E}\rightarrow X$ is a discrete cover with
covering group $\delta _{E}(X)$.
\end{enumerate}
\end{theorem}

\begin{proof}
We have 
\begin{equation*}
\overline{\lambda _{1}}\circ \overline{\lambda _{2}}([\alpha ]_{E})=[\lambda
_{1}\ast \lambda _{2}\ast \alpha ]_{E}=\overline{\lambda _{1}\ast \lambda
_{2}}([\alpha ]_{E})\text{.}
\end{equation*}%
This implies that each $\overline{\lambda }$ is a bijection and therefore is
an element of $H_{X_{E}}$, while the inclusion $\mu :\delta
_{E}(X)\rightarrow H_{X_{E}}$ is a homomorphism. If $\overline{\lambda }%
([\alpha ]_{E})=[\alpha ]_{E}$ for some $\alpha $, then 
\begin{equation*}
\lbrack \ast ]_{E}=[\lambda \ast \alpha ]_{E}\ast \lbrack \alpha
_{E}]^{-1}=[\lambda \ast \alpha \ast \alpha ^{-1}]_{E}=[\lambda ]_{E}\text{.}
\end{equation*}%
This implies both that $\mu $ is injective and the action is free. The
action of $\delta _{E}(X)$ is free. By the very definition, for any
entourage $F\subset E$ and all $\overline{\lambda }\in \delta _{E}(X)$, $%
\overline{\lambda }(F^{\ast })=F^{\ast }$ and therefore the action is
isomorphic. Now suppose that $(x,\overline{\lambda }(x))\in E^{\ast }$ for
some $x\in X_{E}$. Then $\phi _{XE}(x)=\phi _{XE}(\overline{\lambda }(x))$
and since $\phi _{XE}$ is injective on every $E^{\ast }$-ball by Proposition %
\ref{xee}, $\overline{\lambda }(x)=x$. Since we have shown the action is
free, it follows that $\overline{\lambda }$ is the identity and we have
shown that the action is discrete.

For the third part note that by Theorem 11 (see also Remark 13), \cite{P},
we need only check first that $\phi _{XE}$ is bi-uniformly continuous (which
we already know) and second that if $x\in X$ and $y\in X_{E}$ are such that $%
\phi _{XE}(y)=x$, then the preimage $\phi _{XE}^{-1}(x)$ is precisely the
orbit $\delta _{E}(X)(y)$. For any $x\in X$, 
\begin{equation*}
\phi _{XE}^{-1}(x)=\{[\ast =x_{0},...,x_{n-1},x]_{E}\}\text{.}
\end{equation*}%
Now there is some $E$-chain $\alpha =\{\ast =x_{0},...,x_{n-1},x\}$ since $%
\phi _{XE}$ is surjective. Moreover, for any $\lambda $, the endpoint of $%
\overline{\lambda }([\alpha ]_{E})$ is still $x$ and hence $\overline{%
\lambda }([\alpha ]_{E})\in \phi _{XE}^{-1}(x)$. That is, the orbit $\delta
_{E}(X)[\alpha ]_{E}$ is contained in $\phi _{XE}^{-1}(x)$. On the other
hand, if $\beta $ is any other $E$-chain to $x$ and we let $\lambda :=\beta
\ast \alpha ^{-1}$ then $[\lambda ]_{E}\in \delta _{E}(X)$ and $\overline{%
\lambda }([\alpha ]_{E})=[\beta ]_{E}$ and therefore $\phi _{XE}^{-1}(x)$ is
contained in the orbit $\delta _{E}(X)[\alpha ]_{E}$.
\end{proof}

From Proposition \ref{identify} we have:

\begin{corollary}
\label{0530cov}If $X$ is a uniform space then for any entourages $E,F$ in $X$
with $F\subset E$ such that $\phi _{EF}:X_{F}\rightarrow X_{E}$ is
surjective (in particular when $X_{E}$ is chain connected), $\phi _{EF}$ is
a discrete cover with covering group $\delta _{F}(X_{E})$. In particular $%
\phi _{EF}$ is a uniform homeomorphism if and only if $\delta _{F}(X_{E})$
is trivial.
\end{corollary}

\begin{definition}
Let $X$ be a uniform space and $\theta _{EF}:\delta _{F}(X)\rightarrow
\delta _{E}(X)$ denote the restriction of $\phi _{EF}$ to $\delta _{F}(X)$.
The collection $(\theta _{EF},\delta _{F}(X))$ forms an inverse system. We
denote $\underleftarrow{\lim }\delta _{E}(X)$ by $\delta _{1}(X)$ and call
it the deck group of $X$.
\end{definition}

\begin{lemma}
\label{chainhom2}If $X$ is a uniform space and $E$ is an entourage such that 
$X_{E}$ is chain connected then for any entourage $F\subset E$, $\theta
_{EF}:\delta _{F}\rightarrow \delta _{E}$ is surjective.
\end{lemma}

\begin{proof}
According to Lemma \ref{chainhom}, if $\gamma $ is an $E$-loop in $X$, $%
\gamma $ is $E$-homotopic to an $F$-chain $\alpha $. But since $E$-homotopy
preserves endpoints, $\alpha $ must also be an $F$-loop and by construction $%
\theta _{EF}([\alpha ]_{F})=[\alpha ]_{E}=[\gamma ]_{E}$.
\end{proof}

\begin{proposition}
\label{coin}Let $X$ be a uniform space. Then

\begin{enumerate}
\item For any entourage $E$, the group $\delta _{E}(X)$ is discrete with
respect to both the topology of uniform convergence and the topology induced
by the inclusion of $\delta _{E}(X)$ in $X_{E}$.

\item The group $\delta _{1}(X)$ is prodiscrete with respect to the inverse
limit topology, which is the same as the topology induced by the inclusion
in $\widetilde{X}$.
\end{enumerate}
\end{proposition}

\begin{proof}
Since $\delta _{E}(X)$ acts discretely (Theorem \ref{discrete}) it follows
from Corollary 32 of \cite{P} that $\delta _{E}(X)$ is discrete with respect
to the topology of uniform convergence. On the other hand, since $\phi _{XE}$
is injective on $B(\ast ,E^{\ast })$ and $\delta _{E}(X)=\phi
_{XE}^{-1}(\ast )$, $\delta _{E}(X)\cap B(\ast ,E^{\ast })=\ast $, which
shows that $\delta _{E}(X)$ is a discrete subset of $\phi _{XE}$. The second
part now follows from the definitions.
\end{proof}

Before we consider coverable spaces in more detail we need to revisit the
issue of basepoints. Note that the construction of $\widetilde{X}$ is
dependent on the initial choice of basepoint, and in fact $\widetilde{X}$
itself may depend on the choice of the basepoint. For example, if one takes $%
X:=\{0\}\cup \lbrack 1,2]$ with the subspace metric then one can check that $%
\widetilde{X}$ based at $\{0\}$ will consist of a single point, while $%
\widetilde{X}$ based at any point in $[1,2]$ will be $[1,2]$. The following
lemma clears up this issue:

\begin{lemma}
\label{baseok}Suppose that $X$ is a uniform space such that for some choice
of basepoint $\ast $, the projection $\phi :\widetilde{X}\rightarrow X$ is
surjective. Then $X$ is chain connected. Moreover, if $X$ satisfies the
definition of coverable for $\ast $ then

\begin{enumerate}
\item for any two basepoints there is a natural system of uniform
homeomorphisms between the fundamental inverse systems of $X$ with respect
to the two basepoints. In particular the spaces $\widetilde{X}$ constructed
with each basepoint are naturally uniformly homeomorphic.

\item $X$ satisfies the definition of coverable for any choice of basepoint.
\end{enumerate}
\end{lemma}

\begin{proof}
If $\phi :\widetilde{X}\rightarrow X$ is surjective then since $\phi =\phi
_{XE}\circ \phi _{X}$, $\phi _{XE}$ is surjective. It now follows from Lemma %
\ref{lemsur5} that $X$ is chain connected. The first statement now follows
from the observations in Remark \ref{basepoint}. The second statement
follows from the first; essentially the two inverse systems are the same and
so the projections are surjective in one if and only if they are surjective
in the other.
\end{proof}

\begin{theorem}
\label{maincov}Let $X$ be coverable and $\Lambda :=\mathcal{C}(X)$. Then $%
\{X_{E},\phi _{EF}\}_{E\in \Lambda }$ and $\{\delta _{E}(X),\theta
_{EF}\}_{E\in \Lambda }$ comprise an inverse system of discrete covers and $%
\phi :\widetilde{X}\rightarrow X$ is a cover with covering group $\delta
_{1}(X)$.
\end{theorem}

\begin{proof}
Since $X$ is coverable we may take $\Lambda $ for our indexing set. If $%
\lambda $ is an $F$-loop and $\alpha $ is an $F$-chain then 
\begin{equation*}
\phi _{EF}([\lambda ]_{F}([\alpha ]_{F}))=\phi _{EF}([\lambda \ast \alpha
]_{F})=[\lambda \ast \alpha ]_{E}
\end{equation*}%
\begin{equation*}
=[\lambda ]_{E}([\alpha ]_{E})=\theta _{EF}([\lambda ]_{F})(\phi
_{EF}([\alpha ]_{F})\text{.}
\end{equation*}%
That is, the system is compatible in the sense of \cite{P}. Since each $\phi
_{E}$ is surjective, so is each $\phi _{EF}$ and according to Corollary \ref%
{0530cov}, each $\phi _{EF}$ is a discrete cover. The proof is now finished
by \cite{P}, Theorem 44.
\end{proof}

\begin{definition}
When $X$ is coverable we will refer to the projection $\phi :\widetilde{X}%
\rightarrow X$ as the universal covering map of $X$.
\end{definition}

\begin{remark}
We will always take $\delta _{1}(X)$ to have the (prodiscrete) inverse limit
topology, with respect to which $\delta _{1}(X)$ is complete and Hausdorff. $%
\delta _{1}(X)$ is also prodiscrete with respect to the topology of uniform
convergence since $\phi $ is a cover (see \cite{P}). It may be that the
topology of uniform convergence is the the same as the induced topology, but
we have no need for such a statement in this paper.
\end{remark}

\section{Universal and Lifting Properties}

\begin{proposition}
\label{0607uniform}Let $X$ be a uniform space and $\mathcal{U}$ be a
uniformity base. The following are equivalent for a fixed basepoint $\ast $:

\begin{enumerate}
\item $\mathcal{U}$ is universal.

\item $X\times X\in \mathcal{U}$ and for any $F\subset E$ with $E,F\in 
\mathcal{U}$, $\phi _{EF}$ is a uniform homeomorphism.

\item $X\times X\in \mathcal{U}$ and for any $E\in \mathcal{U}$, $\phi _{E}:%
\widetilde{X}\rightarrow X_{E}$ is a uniform homeomorphism.
\end{enumerate}
\end{proposition}

\begin{proof}
Suppose that $\mathcal{U}$ is universal and $F\subset E$ are entourages in $%
\mathcal{U}$. Then $\phi _{EF}:X_{F}\rightarrow X_{E}$ satisfies $\phi
_{EF}=\phi _{XE}^{-1}\circ \phi _{XF}$ and is therefore a uniform
homeomorphism. Since $X\times X$ lies in any universal base, 1)$\Rightarrow $%
2). If 2) holds then we may use $\mathcal{U}$ as our indexing set for the
fundamental inverse system of $X$ and since each of the bonding maps is a
uniform homeomorphism the inverse system is in fact trivial. Therefore each $%
\phi _{E}:\widetilde{X}\rightarrow X_{E}$ is a uniform homeomorphism and 3)
follows. If 3) is given then note that for any $E\in \mathcal{U}$ we have
that $\phi _{XE}=\phi _{X}\circ \phi _{E}^{-1}$ is a uniform homeomorphism
and we see that $\mathcal{U}$ is universal.
\end{proof}

Definition \ref{universaldef} implicitly involves the choice of a basepoint
for the construction of the spaces $X_{E}$. However, as the next Corollary
(the proof of which is immediate from the Proposition \ref{0607uniform} and
Lemma \ref{baseok}) shows that the definition is independent of basepoint.

\begin{corollary}
\label{uic}Any universal space is coverable, hence chain connected. In
particular, if $X$ is universal with respect to one basepoint then $X$ is
universal with respect to any basepoint.
\end{corollary}

\begin{corollary}
\label{0518u}If $X$ is universal and $E$ is an entourage such that $X_{E}$
is chain connected then $\phi _{XE}:X_{E}\rightarrow X$ is a uniform
homeomorphism.
\end{corollary}

\begin{proof}
Let $F\subset E$ be an entourage in the universal base. Since $X_{E}$ is
chain connected, $\phi _{EF}:X_{F}\rightarrow X_{E}$ is surjective. Since $%
\phi _{XF}:X_{F}\rightarrow X$ is bijective and $\phi _{XF}=\phi _{XE}\circ
\phi _{EF}$, $\phi _{EF}$ must also be injective. Since $X_{E}$ is chain
connected, Lemma \ref{lemsur5} implies that $\phi _{EF}$ is surjective,
hence a uniform homeomorphism because $\phi _{EF}$ is bi-uniformly
continuous. Therefore $\phi _{XE}=\phi _{XF}\circ \phi _{EF}^{-1}$ is a
uniform homeomorphism.
\end{proof}

\begin{theorem}
\label{xic}If $X$ is coverable then $\widetilde{X}$ is universal.
\end{theorem}

\begin{proof}
Applying Proposition \ref{coverchain} and Lemma \ref{41} we have that $%
\widetilde{X}$ is chain connected. A basis element for the uniformity of $%
\widetilde{X}$ is of the form $\phi _{E}^{-1}(F^{\ast })$, where $F\subset E$
are entourages in $X$. Let $D$ be a covering entourage contained in $E\cap F$%
. Then $\phi _{D}$ is surjective and 
\begin{equation*}
\phi _{E}(\phi _{D}^{-1}(D^{\ast }))=\phi _{ED}(\phi _{D}(\phi
_{D}^{-1}(D^{\ast }))
\end{equation*}%
\begin{equation*}
=\phi _{ED}(D^{\ast })=D^{\ast }\subset F^{\ast }\subset X_{E}\times X_{E}%
\text{.}
\end{equation*}%
In other words, $\phi _{D}^{-1}(D^{\ast })\subset \phi _{E}^{-1}(F^{\ast })$
and we may take for our basis elements of $\widetilde{X}$ entourages of the
form $G:=\phi _{D}^{-1}(D^{\ast })$. We will now apply Corollary \ref{67},
taking $f$ to be the surjective map $\phi _{D}:\widetilde{X}\rightarrow
X_{D} $. We have the following diagram: 
\begin{equation*}
\begin{array}{lll}
\widetilde{X}_{G} &  & (X_{D})_{D^{\ast }} \\ 
\downarrow ^{\phi _{\widetilde{X}G}} &  & \downarrow ^{\phi _{X_{D}D^{\ast
}}} \\ 
\widetilde{X} & \overset{\phi _{D}}{\longrightarrow } & X_{D}%
\end{array}%
\end{equation*}%
According to Proposition \ref{identify} the function $\phi _{X_{D}D^{\ast
}}:(X_{D})_{D^{\ast }}\rightarrow X_{D}$ is a uniform homeomorphism and we
may conclude that $\phi _{\widetilde{X}G}$ is also a uniform homeomorphism.
The collection of all such $G$ is therefore a universal base for $\widetilde{%
X}$.
\end{proof}

\begin{corollary}
\label{uit}The following are equivalent for a coverable space $X$.

\begin{enumerate}
\item $X$ is universal.

\item For every entourage $E$ such that $X_{E}$ is chain connected, $\delta
_{E}(X)$ is trivial (i.e., every $E$-loop based at $\ast $ is $E$-homotopic
to the trivial loop).

\item $\delta _{1}(X)$ is trivial.
\end{enumerate}
\end{corollary}

\begin{proof}
The implication (1)$\Rightarrow $(2) follows from Corollary \ref{0518u}. If
(2) holds then for any coverable entourage $E$, $\phi _{E}:\widetilde{X}%
\rightarrow X_{E}$ is a uniformly contionuous surjection and hence $X_{E}$
is chain connected. Therefore each $\delta _{E}(X)$ is trivial and the
inverse limit $\delta _{1}(X)$ is trivial. Finally, if $\delta _{1}(X)$ is
trivial then the universal covering map is a uniform homeomorphism and since 
$\widetilde{X}$ is universal, so is $X$.
\end{proof}

\begin{example}
We will see later (Theorem \ref{0627simply}) that $\mathbb{R}$ is universal.
Consider the entourage $E$ consisting of all $(x,y)$ such that $x-y\in W$,
where $W:=(-1,1)\cup (2,4)\cup (-4,-2)$. It is not hard to see that $%
\{0,3,0\}$ is an $E$-loop based at $0$ that is not $E$-homotopic to the
trivial loop. At the same time it is true that $\mathbb{R}_{E}$ is uniformly
homeomorphic to $\mathbb{R\times Z}$, where $\mathbb{Z}$ has the discrete
uniformity, and hence $\mathbb{R}_{E}$ is not chain connected (see Example
48 in \cite{BPTG} for more details). This shows that one cannot expect $%
\delta _{E}(X)$ to be trivial for every choice of $E$ when $X$ is universal.
\end{example}

\begin{proposition}
If $f:X\rightarrow Y$ is a uniformly continuous bijection between coverable
spaces $X$ and $Y$ and $Y$ is universal then $X$ is universal.
\end{proposition}

\begin{proof}
Suppose $Y$ is universal and let $F$ be an entourage in the universal base
for $Y$. Then $\phi _{YF}:Y_{F}\rightarrow Y$ is a uniform homeomorphism and
by Corollary \ref{67} $\phi _{XE}:X_{E}\rightarrow X$ is a uniform
homeomorphism, where $E:=f^{-1}(F)$.
\end{proof}

\begin{proposition}
\label{disccover}Let $f:X\rightarrow X/G=Y$ be a discrete cover, where $X$
is chain connected. For any sufficiently small $G$-invariant entourage $E$
and $F:=f(E)$, the function $f_{EF}:X_{E}\rightarrow Y_{F}$ is a uniform
homeomorphism.
\end{proposition}

\begin{proof}
Since $f$ is a discrete cover, there exists an entourage $D$ such that if $%
g\in G$ satisfies $(g(x),x)\in D$ for some $x\in X$ then $g=e$. Suppose that 
$E$ is any invariant entourage such that $E^{3}\subset D$. Since $f$ is a
quotient mapping, $F=f(E)$ is an entourage in $Y$. We will first show that $%
f_{EF}$ is injective. By Lemma \ref{lemsur2} we may equivalently prove that
if $\gamma =\{x_{0},...,x_{n}\}$ is an $E$-chain in $X$ such that $%
x_{n}=g(x_{0})$ for some $g\in G$ and the $F$-loop $f(\gamma )$ is $F$%
-homotopic to the trivial loop $\{f(x_{0})\}$ then $g=e$ (so $\gamma $ is an 
$E$-loop) and $\gamma $ is $E$-homotopic to the trivial loop $\{x_{0}\}$.
Let $f(\gamma )=\{y_{0},...,y_{n}\}$. We will prove the statement by
induction on the minimal length $m$ of an $F$-homotopy between $f(\gamma )$
and $\{y_{0}\}$. If $m=0$ then $f(\gamma )$ is already trivial and so is $%
\gamma $, and the proof is finished. Suppose we have proved it for some $%
m-1\geq 0$, and there is some $F$-homotopy of $f(\gamma )$ to $\{y_{0}\}$ of
length $m$. Suppose that the first step in the $F$-homotopy is to add a
point: $f(\gamma )$ is $F$-related to $\{y_{0},...,y_{k},y,y_{k+1},...,y_{n}%
\}$. That is, $(y_{k},y),(y,y_{k+1})\in F=f(E)$. Now there exists $(a,b)\in
E $ such that $f(a)=y_{k}$ and $f(b)=y$, and therefore some $g\in G$ such
that $g(a)=x_{k}$. Since $E$ is invariant, we have that if $w:=g(b)$ then $%
(w,x_{k})\in E$ and $f(w)=y$. Likewise there is some $w^{\prime }\in X$ such
that $f(w^{\prime })=y$ and $(x_{k+1},w^{\prime })\in E$. Since $%
(x_{k},x_{k+1})\in E$, $(w,x_{k})\in E$, and $(x_{k+1},w^{\prime })\in E$,
it follows that $(w,w^{\prime })\in E^{3}\subset D$. But since $w,w^{\prime
}\in f^{-1}(y)$, $w=k(w^{\prime })$ for some $k\in G$ and by choice of $D$,
it must be that $k=e$ and $w=w^{\prime }$. The inductive hypothesis now
finishes the proof. Now suppose that the first step in the homotopy is to
remove a point: $f(\gamma )$ is $F$-related to $%
\{y_{0},...,y_{k-1},y_{k+1},...,y_{n}\}$; that is, $\left(
y_{k-1},y_{k+1}\right) \in F$. As in the preceding argument there is some $%
g\in G$ such that $(x_{k-1},g(x_{k+1}))\in E$. But $%
(x_{k-1},x_{k}),(x_{k},x_{k+1})\in E$, so $(g(x_{k+1}),x_{k+1})\in
E^{3}\subset D$, and we conclude that $g(x_{k+1}))=x_{k+1}$. Therefore $%
(x_{k-1},x_{k+1})\in E$ and we may again apply the inductive hypothesis.

To see why $f_{EF}$ is surjective, let $\{\ast =y_{0},...,y_{n}\}$ be an $F$%
-chain in $Y$. We will prove the statement by induction in $n$. If $n=0$
then the proof is obvious since $\ast =f(\ast )$ by assumption. Suppose we
have proved it for $n-1\geq 0$. Then we can find some $E$-chain $\{\ast
=x_{0},...,x_{n-1}\}$ such that $y_{i}=f(x_{i})$ for all $i$. Since $%
(y_{n-1},y_{n})\in F=f(E)$, and $f(x_{n-1})=y_{n-1}$, we may again use the
invariance of $E$ to see that there exists some $x_{n}\in X$ such that $%
f(x_{n})=y_{n}$ and $(x_{n-1},x_{n})\in E$. Then $\{x_{0},...,x_{n}\}$ is an 
$E$-chain such that $f(\{x_{0},...,x_{n}\})=\{y_{0},...,y_{n}\}$.

To finish the proof of the proposition, let $D^{\ast }$ be an entourage in $%
X_{E}$, where $D\subset E$ is an invariant entourage in $X$ ($D$ exists
since $G$ acts isomorphically). Since $X$ is chain connected, $\phi
_{XE}(D^{\ast })=D$ by Proposition \ref{xee}, and $f(D)=K\subset F$ is an
entourage in $X/G$. The proof will be complete if we show that $K^{\ast
}\subset f_{EF}(D^{\ast })$, which makes $f_{EF}(D^{\ast })$ an entourage.
Let $([\alpha ]_{F},[\beta ]_{F})\in K^{\ast }$ which means 
\begin{equation*}
([\alpha ]_{F},[\beta ]_{F})=([\ast =y_{0},...,y_{n},y]_{F},[\ast
=y_{0},...,y_{n},z]_{F})
\end{equation*}%
with $(y,z)\in K=f(D)$. Using the invariance of $E$ and $D$ and proceeding
inductively as we have done above, we can find $E$-chains $\alpha ^{\prime
}=\left\{ \ast =x_{0},...,x_{n},y^{\prime }\right\} $ and $\beta ^{\prime
}=\left\{ \ast =x_{0},...,x_{n},z^{\prime }\right\} $ such that $f(\alpha
^{\prime })=\alpha $ and $f(\beta ^{\prime })=\beta $, and $(x_{n},y^{\prime
}),(x_{n},z^{\prime })\in D$. By definition, $([\alpha ^{\prime
}]_{E},[\beta ^{\prime }]_{E})\in D^{\ast }$ and $\left( f_{EF}([\alpha
^{\prime }]_{E}),f_{EF}([\alpha ^{\prime }]_{E})\right) =([\alpha
]_{F},[\beta ]_{F})$.
\end{proof}

\begin{remark}
\label{base2}We give one final comment about basepoints and lifts of
functions following Lemma \ref{baseok}. The lifting theorems below are true
for any choice of basepoints such that the functions involved are
basepoint-preserving. For example, in the proposition below we may start
with a basepoint $\ast $ in $X$, choose $\ast =f(\ast )$ in $Y$, use any
basepoint $\ast ^{\prime }$ to construct $\widetilde{Y}$ and then choose
another basepoint $\ast $ in $\widetilde{Y}$ so that $\phi (\ast )=\ast $.
\end{remark}

\begin{proposition}
\label{liftlem}Let $X$ be universal, $Y$ be uniform and $f:X\rightarrow Y$
be uniformly continuous. Then

\begin{enumerate}
\item For any entourage $E$ in $Y$ there is a unique uniformly continuous
function $f_{E}:X\rightarrow Y_{E}$ such that $\phi _{YE}\circ f_{E}=f$ and $%
f_{E}(\ast )=\ast $.

\item There is a unique uniformly continuous function $f_{L}:X\rightarrow 
\widetilde{Y}$ such that $f_{L}(\ast )=\ast $ and $\phi \circ f_{L}=f$,
where $\phi :\widetilde{Y}\rightarrow Y$ is the projection.
\end{enumerate}
\end{proposition}

\begin{proof}
Define $f_{E}:X\rightarrow Y_{E}$ as follows. Let $F$ be an entourage in the
universal base of $X$ such that $f(F)\subset E$ and $k_{E}:X_{F}\rightarrow
Y_{E}$ be the unique uniformly continuous function given by Theorem \ref{60}
($\phi _{XF}:X_{F}\rightarrow X$ is a uniform homeomorphism and therefore $%
X_{F}$ is chain connected). Define $f_{E}:=k_{E}\circ \phi _{XF}^{-1}$. If $%
g $ were any such function then it follows from the uniqueness of $k_{E}$
that $k_{E}=g\circ \phi _{XF}$ and hence that $g=f_{E}$.

Note that by uniqueness, if $E\subset F$ are entourages in $Y$ we have that $%
\phi _{FE}\circ f_{E}=f_{F}$ and by the universal property of the inverse
limit there is a unique function $f_{L}:X\rightarrow \widetilde{Y}$ such
that $\phi _{E}\circ f_{L}=f_{E}$ for every entourage $E$ and $f_{L}(\ast
)=\ast $. Suppose that $f^{\prime }:X\rightarrow \widetilde{Y}$ is uniformly
continuous such that $f^{\prime }(\ast )=\ast $ and $\phi \circ f^{\prime
}=f $. Note that for any entourage $E$ in $Y$ we have 
\begin{equation*}
\phi _{YE}\circ \left( \phi _{E}\circ f^{\prime }\right) =\phi \circ
f^{\prime }=f
\end{equation*}%
and therefore by uniqueness in Part (1), $\phi _{E}\circ f^{\prime }=f_{E}$.
Since $f^{\prime }$ is also induced by the functions $\phi _{E}\circ
f^{\prime }$, $f=f^{\prime }$.
\end{proof}

\begin{notation}
The functions $f_{E}$ and $f_{L}$ will both be referred to as
\textquotedblleft lifts\textquotedblright\ of $f$.
\end{notation}

\begin{theorem}
\label{76}Let $X$ and $Y$ be uniform, $f:X\rightarrow Y$ be a cover, $Z$ be
universal, and $g:Z\rightarrow Y$ be uniformly continuous. Then there exists
a unique uniformly continuous function $h:Z\rightarrow X$ such that $f\circ
h=g$ and $h(\ast )=\ast $.
\end{theorem}

\begin{proof}
Suppose first that $f$ is a discrete cover. According to Proposition \ref%
{disccover} there is an entourage $E$ in $X$ such that if $F:=f(E)$ then $%
f_{EF}:X_{E}\rightarrow Y_{F}$ is a uniform homeomorphism and $f\circ \phi
_{XE}=\phi _{YF}\circ f_{EF}$. Define $h:Z\rightarrow X$ by $\phi _{XE}\circ
f_{EF}^{-1}\circ g_{F}$ (where $g_{F}$ is the lift of $g$ given by
Proposition \ref{liftlem}). Then $h(\ast )=\ast $ and 
\begin{equation*}
f\circ h=f\circ \phi _{XE}\circ f_{EF}^{-1}\circ g_{F}
\end{equation*}%
\begin{equation*}
=\phi _{YF}\circ f_{EF}\circ f_{EF}^{-1}\circ g_{F}=\phi _{YF}\circ g_{F}=g%
\text{.}
\end{equation*}%
To prove uniqueness, suppose that $h^{\prime }$ is any such function.
Consider the lift $h_{E}^{\prime }:Z\rightarrow X_{E}$. We have 
\begin{equation*}
\phi _{YF}\circ \left( f_{EF}\circ h_{E}^{\prime }\right) =f\circ \phi
_{XE}\circ h_{E}^{\prime }=f\circ h^{\prime }=g\text{.}
\end{equation*}%
By uniqueness of lifts (Proposition \ref{liftlem}), $f_{EF}\circ
h_{E}^{\prime }=g_{F}$. But then $h_{E}^{\prime }=f_{EF}^{-1}\circ g_{F}$
and 
\begin{equation*}
h^{\prime }=\phi _{XE}\circ h_{E}^{\prime }=\phi _{XE}\circ f_{EF}^{-1}\circ
g_{F}=h\text{.}
\end{equation*}

Now suppose that $f$ is an arbitrary cover. By Theorem 48 in \cite{P} there
exists an inverse system $\{X_{\alpha },f_{\alpha \beta }\}$ such that $%
f_{\alpha \beta }:X_{\beta }\rightarrow X_{\alpha }$ is a discrete cover and 
$X=\underleftarrow{\lim }X_{\alpha }$ and $Y=X_{_{1}}$ for some minimal
element $1$. According to what we proved above, for each $\alpha $ there is
a unique uniformly continuous function $h_{\alpha }:Z\rightarrow X_{\alpha }$
such that $h_{\alpha }(\ast )=\ast $ and $f_{1\alpha }\circ h_{\alpha }=g$.
If $\alpha \leq \beta $ we have that 
\begin{equation*}
f_{1\alpha }\circ \left( f_{\alpha \beta }\circ h_{\beta }\right) =f_{1\beta
}\circ h_{\beta }=g
\end{equation*}%
and by uniqueness $f_{\alpha \beta }\circ h_{\beta }=h_{\alpha }$. By the
universal property of inverse limits there is a unique uniformly continuous
function $h:Z\rightarrow X=\underleftarrow{\lim }X_{\alpha }$ such that for
all $\alpha $, $f_{\alpha }\circ h=h_{\alpha }$ and $h(\ast )=\ast $. Now
suppose $h^{\prime }:Z\rightarrow X$ is any uniformly continuous function
such that $h^{\prime }(\ast )=\ast $ and $f\circ h^{\prime }=g$. Define $%
h_{\alpha }^{\prime }:=f_{\alpha }\circ h^{\prime }$. We have that $%
h_{\alpha }^{\prime }(\ast )=\ast $ and 
\begin{equation*}
f_{1\alpha }\circ h_{\alpha }^{\prime }=f_{1\alpha }\circ f_{\alpha }\circ
h^{\prime }=f\circ h^{\prime }=g\text{.}
\end{equation*}%
By the uniqueness of $h_{\alpha }$, $h_{\alpha }=h_{\alpha }^{\prime }$ and
therefore $h=h^{\prime }$.
\end{proof}

If $Y$ is coverable then by Theorem \ref{xic}, $\widetilde{Y}$ is universal
and we obtain:

\begin{corollary}
\label{partialstep}Let $f:X\rightarrow Y$ be a cover where $X$ is uniform
and $Y$ is coverable, and $\phi :\widetilde{Y}\rightarrow Y$ be the
projection. Then there exists a unique uniformly continuous function $f_{B}:%
\widetilde{Y}\rightarrow X$ such that $f_{B}(\ast )=\ast $ and $f\circ
f_{B}=\phi $.
\end{corollary}

\begin{corollary}
\label{quick}If $f:X\rightarrow Y$ is a bi-uniformly continuous surjection
where $X$ is universal and $Y$ is uniform then $Y$ is coverable.
\end{corollary}

\begin{proof}
Let $\mathcal{U}$ be a universal base for $X$. Since $f$ is a bi-uniformly
continuous surjection, the set $\mathcal{B}$ of all $f(F)$ such that $F\in 
\mathcal{U}$ is a base for the uniformity of $Y$. We may index the
fundamental system for $Y$ using $\mathcal{B}$, and for each $E=f(F)$ in $%
\mathcal{B}$ we may use $F$ in the construction of the lift $f_{E}$. Since $%
f_{E}=\phi _{E}\circ f_{L}$ we will be finished by the definition of
coverable if we can show that each $f_{E}$ is surjective. Let $\beta =\{\ast
=y_{0},...,y_{n}\}$ be an $E$-chain in $Y$. We will show by induction on $n$
that there is an $F$-chain $\alpha =\{\ast =x_{0},...,x_{n}\}$ such that $%
f(\alpha )=\beta $ and, by definition of $f_{E}$, this completes the proof.
For $n=0$ the proof is obvious. Now suppose we have an $F$-chain $\gamma
=\{\ast =x_{0},...,x_{n-1}\}$ such that $f(\gamma )=\{\ast
=y_{0},...,y_{n-1}\}$. Since $f(F)=E$ there is some ordered pair $%
(x_{n-1},x_{n})\in F$ such that $f((x_{n-1},x_{n}))=(y_{n-1},y_{n})$. The $F$%
-chain $\alpha =\{\ast =x_{0},...,x_{n-1}\}$ now satisfies $f(\alpha )=\beta 
$.
\end{proof}

\begin{corollary}
\label{bi}If $X$ is coverable and $f:X\rightarrow Y$ is a bi-uniformly
continuous surjection then $Y$ is coverable. In particular, any quotient by
an equiuniform action, hence by an isomorphic action, on a coverable uniform
space is coverable.
\end{corollary}

\begin{proof}
Since the universal cover $\phi :\widetilde{X}\rightarrow X$ is a
bi-uniformly continuous surjection, so is $f\circ \phi :\widetilde{X}%
\rightarrow Y$. The proof is finished by Corollary \ref{quick} and Theorem %
\ref{xic}.
\end{proof}

\begin{theorem}
\label{73}Let $X,Y$ be coverable spaces, $f:X\rightarrow Y$ be uniformly
continuous, and $\phi :\widetilde{X}\rightarrow X$ and $\psi :\widetilde{Y}%
\rightarrow Y$ be the projections. Then there is a unique uniformly
continuous function $\widetilde{f}:\widetilde{X}\rightarrow \widetilde{Y}$,
such that $\widetilde{f}(\ast )=\ast $ and $f\circ \phi =\psi \circ 
\widetilde{f}$. Moreover,

\begin{enumerate}
\item For any $x\in \widetilde{X}$ and $g\in \delta _{1}(X)$, $\widetilde{f}%
(g(x))=\widetilde{f}(g)(\widetilde{f}(x))$.

\item The restriction $f_{\ast }$ of $\widetilde{f}$ to $\delta _{1}(X)$ is
a homomorphism into $\delta _{1}(Y)$.

\item If $f$ is a discrete cover then $\widetilde{f}$ is a uniform
homeomorphism.

\item If $Z$ is uniform and $g:Y\rightarrow Z$ is uniformly continuous then $%
\widetilde{g\circ f}=\widetilde{g}\circ \widetilde{f}$ (and in particular $%
\left( g\circ f\right) _{\ast }=g_{\ast }\circ f_{\ast }$).
\end{enumerate}
\end{theorem}

\begin{proof}
For the main statement, define $\widetilde{f}:=\left( f\circ \phi \right)
_{L}$. If $g:\widetilde{X}\rightarrow \widetilde{Y}$ is a uniformly
continuous function with $\psi \circ g=f\circ \phi $ then by definition $g$
is a lift of $f\circ \phi $ and so $g=\widetilde{f}$.

For Part (1), note that if $x\in \delta _{1}(X)$ then $\phi (x)=\ast $ and 
\begin{equation*}
\psi \circ \widetilde{f}(x)=f\circ \phi (x)=f(\ast )=\ast \text{.}
\end{equation*}%
Therefore $\widetilde{f}(\delta _{1}(X))\subset \delta _{1}(Y)$. For any
entourage $F$ in $Y$ let $E$ be an entourage in $X$ such that $f(E)\subset F$%
. Now $\psi _{F}\circ \widetilde{f}$ and $f_{EF}\circ \phi _{E}$ are both
lifts of $f\circ \phi $ to $Y_{F}$ and therefore $\psi _{F}\circ \widetilde{%
f=}f_{EF}\circ \phi _{E}$ (see Proposition \ref{liftlem}). Now let $g\in
\delta _{1}(X)$ and $x\in \widetilde{X}$ with $\phi _{E}(g)=\left[ \gamma
_{E}\right] _{E}$ and $\phi _{E}(x)=\left[ \alpha _{E}\right] _{E}$. By
definition of the action of the inverse limit group $G$ on the inverse limit
space $\widetilde{X}$,%
\begin{equation*}
\psi _{F}(\widetilde{f}(g(x)))=f_{EF}\circ \phi _{E}(g(x))=f_{EF}\left( %
\left[ \gamma _{E}\ast \alpha _{E}\right] _{E}\right) \text{.}
\end{equation*}%
The latter quantity, by definition of $f_{EF}$, is equal to%
\begin{equation*}
\left[ f(\gamma _{E}\ast \alpha _{E})\right] _{F}=\left[ f(\gamma _{E})\ast
f(\alpha _{E})\right] _{F}=\left[ f(\gamma _{E})\right] _{F}\ast \left[
f(\alpha _{E})\right] _{F}\text{.}
\end{equation*}%
As a special case when the $g$ is the identity, 
\begin{equation*}
\psi _{F}(\widetilde{f}(x))=\left( [f(\alpha _{E})]_{F}\right) \text{.}
\end{equation*}%
If $x=\ast $ then 
\begin{equation*}
\psi _{F}(\widetilde{f}(g(\ast )))=\left[ f(\gamma _{E})\right] _{F}\text{.}
\end{equation*}%
Combining these we obtain 
\begin{equation*}
\widetilde{f}(g(x))=\widetilde{f}(g)(\widetilde{f}(x))\text{.}
\end{equation*}%
When $x\in \delta _{1}(X)$, we have%
\begin{equation*}
\widetilde{f}(gx)=\widetilde{f}(g)\widetilde{f}(x)\text{,}
\end{equation*}%
which gives the second statement.

If $f:X\rightarrow X/G=Y$ is a discrete cover then $f$ is bi-uniformly
continuous and we may index the fundamental system of $Y$ using entourages
of the form $f(E)$ where $E$ is an invariant entourage in $X$. Then by
uniqueness, $\widetilde{f}$ is induced by the functions $f_{EF}$.
Proposition \ref{disccover} implies that the functions $f_{EF}$ are all
uniform homeomorphisms and hence $\widetilde{f}$ must be a uniform
homeomorphism.

The last part follows from uniqueness of $\widetilde{g\circ f}$.
\end{proof}

\begin{theorem}
\label{101}If $X$ and $Y$ are coverable spaces and $f:X\rightarrow Y$ is a
cover then $\widetilde{f}:\widetilde{X}\rightarrow \widetilde{Y}$ is a
uniform homeomorphism. If $X$ is coverable then $f_{B}:\widetilde{Y}%
\rightarrow X$ is a cover with covering group $f_{\ast }(\delta
_{1}(X))\subset \delta _{1}(Y)$.
\end{theorem}

\begin{proof}
Let $\phi :\widetilde{X}\rightarrow X$ and $\psi :\widetilde{Y}\rightarrow Y$
be the projections. We will show that the lift $(f_{B})_{L}:\widetilde{Y}%
\rightarrow \widetilde{X}$ given by Corollary \ref{partialstep} and
Proposition \ref{liftlem} is an inverse to $\widetilde{f}$. First note that 
\begin{equation*}
\psi \circ \left( \widetilde{f}\circ (f_{B})_{L}\right) =f\circ \phi \circ
(f_{B})_{L}=f\circ f_{B}=\psi \text{.}
\end{equation*}%
That is, $\widetilde{f}\circ (f_{B})_{L}$ is the unique lift of the identity
on $Y$ and hence must be the identity on $\widetilde{Y}$.

Now according to Theorem \ref{76}, $\phi $ is the unique lift of the
function $f\circ \phi $; that is, the unique uniformly continuous function $%
\eta :\widetilde{X}\rightarrow X$ such that such that $\eta (\ast )=\ast $
and $f\circ \eta =f\circ \phi $. But we also have that 
\begin{equation*}
f\circ \left( f_{B}\circ \widetilde{f}\right) =\psi \circ \widetilde{f}%
=f\circ \phi \text{.}
\end{equation*}%
In other words, $f_{B}\circ \widetilde{f}=\phi $. We now have 
\begin{equation*}
\phi \circ ((f_{B})_{L}\circ \widetilde{f})=f_{B}\circ \widetilde{f}=\phi
\end{equation*}%
and $\left( f_{B}\right) _{L}\circ \widetilde{f}$ must be the unique lift of
the identity on $X$, hence the identity on $\widetilde{X}$. We have shown
that both $\left( f_{B}\right) _{L}$ and $\widetilde{f}$ are uniform
homeomorphisms and inverses of one another.

For the second part note that since $X$ is coverable, $\phi $ is
bi-uniformly continuous and since $f_{B}=\phi \circ \left( f_{B}\right)
_{L}=\phi \circ \widetilde{f}^{-1}$, $f_{B}$ is bi-uniformly continuous. Now 
$f_{B}(x)=f_{B}(y)$ for $x,y\in \widetilde{Y}$ if and only $\phi \circ 
\widetilde{f}^{-1}(x)=\phi \circ \widetilde{f}^{-1}(y)$, or equivalently,
letting $w:=\widetilde{f}^{-1}(x)$ and $z:=\widetilde{f}^{-1}(y)$ there is
some $g\in \delta _{1}(X)$ such that $g(w)=z$. But this is equivalent to 
\begin{equation*}
f_{\ast }(g)(x)=f_{\ast }(g)(\widetilde{f}(w))=\widetilde{f}(g(w))=%
\widetilde{f}(z)=y
\end{equation*}%
(the second equality comes from Theorem \ref{73} (1)). That is, the orbits
of $f_{\ast }(\delta _{1}(X))$ are precisely the preimages of points with
respect to $f_{B}$. According to \cite{P} this means that $f_{B}$ is the
quotient with respect to the action of $f_{\ast }(\delta _{1}(X))$. Finally,
since $f_{\ast }(\delta _{1}(X))\subset \delta _{1}(Y)$, which acts
prodiscretely and isomorphically, so does $f_{\ast }(\delta _{1}(X))$ and $%
f_{B}$ is a cover with covering group $f_{\ast }(\delta _{1}(X))$ (see
Remark 13 in \cite{P}).
\end{proof}

Since a universal space is uniformly homeomorphic to its own universal cover
by Proposition \ref{0607uniform}, we obtain:

\begin{corollary}
The universal cover of a coverable space is unique up to uniform
homeomorphism. More precisely, if $X$ is coverable, $Y$ is universal, $%
f:Y\rightarrow X$ is a cover and $\phi :\widetilde{X}\rightarrow X$ is the
universal cover then $f_{L}:Y\rightarrow \widetilde{X}$ is a uniform
homeomorphism.
\end{corollary}

\section{Traditional Topological Properties}

We will say that a topological space $X$ is simply connected if every loop
in $X$ is null-homotopic, regardless of whether $X$ is pathwise connected.
Recall that $X$ is called locally connected (resp. locally pathwise
connected) if for every $x\in X$ and open set $U$ containing $x$ there is a
connected (resp. pathwise connected) open set $V$ with $x\in V\subset U$. $X$
is semilocally simply connected if each $x\in X$ is contained in an open set 
$U$ such that every loop in $U$ based at $x$ is null-homotopic in $X$ (see 
\cite{M}). In a uniform space it is natural to consider the situation when
these local conditions are true uniformly.

\begin{definition}
\label{unidef}A uniform space $X$ is called uniformly locally connected
(resp. uniformly locally pathwise connected) if for each entourage $E$ there
is an entourage $F\subset E$ such that all $F$-balls are open and connected
(resp. open and pathwise connected). $X$ is called uniformly semilocally
simply connected if there exists an entourage $E$ such that any loop in $%
B(x,E)$ based at $x\in X$ is null-homotopic in $X$.
\end{definition}

Note that given any $E$ as in the above definition, any entourage $F\subset
E $ has the same property.

\begin{proposition}
A uniform space $X$ is uniformly locally connected (resp. uniformly locally
pathwise connected) if and only if for every entourage $E$ there exists an
entourage $F$ with open balls such that for every $x\in X$ there exists some
open set $U_{x}$ such that $B(x,F)\subset U_{x}\subset B(x,E)$ and $U_{x}$
is connected (resp. $U$ is pathwise connected).
\end{proposition}

\begin{proof}
Necessity is trivial; take an entourage $F$ with open balls and let $%
U_{x}:=B(x,F)$. To prove the converse, let $E$ be any entourage in $X$ and $%
K $ be an entourage such that $K^{2}\subset E$. Let $D$ be an entourage such
that for each $x\in X$ there exists an open set $U_{x}$ that is connected
(resp. pathwise connected), such that $B(x,D)\subset U_{x}\subset B(x,K)$.
Define $F\subset X\times X$ by 
\begin{equation*}
F:=\{(x,y):x,y\in U_{w}\text{ for some }w\}\text{.}
\end{equation*}%
$F$ is symmetric by definition and certainly contains the diagonal. Now 
\begin{equation*}
B(x,F)=\{y:y\in U_{w}\text{ for some }U_{w}\text{ containing }x\}
\end{equation*}%
which, being the union of connected (resp. pathwise connected) open sets $%
U_{w}$ containing $x$, is open and connected (resp. pathwise connected). $F$
is an entourage since by definition $D\subset F$. Finally, if $(x,y)\in F$
then $y\in B(x,F)$ and $y\in U_{w}$ for some $U_{w}$ containing $x$. Since $%
U_{w}\subset B(w,K)$, $(x,y)\in K^{2}\subset E$.
\end{proof}

\begin{proposition}
\label{428unif}Let $X$ be a locally connected (resp. locally pathwise
connected) compact topological space. Then $X$, with the unique uniformity
compatible with its topology, is uniformly locally connected (resp.
uniformly locally pathwise connected). If $X$ is both locally pathwise
connected and semilocally simply connected then $X$ is uniformly semilocally
simply connected.
\end{proposition}

\begin{proof}
Recall that the unique compatible uniformity on $X$ has as a basis all
symmetric open subsets of $X\times X$ that contain the diagonal. Let $E$ be
any entourage in $X$ and for any $x\in X$ let $U_{x}$ be an open and
connected (resp. pathwise connected) neighborhood of $x$ such that $%
U_{x}\times U_{x}\subset E$. Let $F\subset E$ be the union of all sets $%
U_{x}\times U_{x}$. Then $F$ is open and symmetric, hence an entourage. If $%
y\in B(x,F)$ then $(x,y)\in U_{z}\times U_{z}$ for some $z$, and therefore $%
x $ and $y$ both lie in $U_{z}$. That is, $B(x,F)$ is the union of all sets $%
U_{z}$ such that $x\in U_{z}$. Since each $U_{z}$ is connected (resp.
pathwise connected), so is $B(x,F)$.

If $X$ is both locally pathwise connected and semilocally simply connected
then by what we proved above, for every $x\in X$ there is an arbitrarily
small entourage $E_{x}$ such that the ball $B(x,E_{x})$ is open and pathwise
connected. By choosing $E_{x}$ small enough we may assume that every loop in 
$B(x,E_{x})$ based at $x$ is null-homotopic in $X$. For each $x$ let $F_{x}$
be an entourage with open balls such that $F_{x}^{2}\subset E_{x}$. Let $%
\{B(x_{i},F_{x_{i}})\}_{i=1}^{k}$ be a finite open cover of $X$, and let $%
F:=\bigcap\limits_{i=1}^{k}F_{x_{i}}$. For any $x\in X$ there is some $x_{i}$
such that $(x,x_{i})\in F_{x_{i}}$. If $y\in B(x,F)$ then $(x,y)\in F\subset
F_{x_{i}}$ and $(y,x_{i})\in F_{x_{i}}^{2}\subset E_{x_{i}}$. That is, $%
B(x,F)\subset B(x_{i},E_{x_{i}})$. Given any loop $\gamma $ in $B(x,F)$
based at $x$, join it to $x_{i}$ by a path in $B(x_{i},E_{x_{i}})$ from $x$
to $x_{i}$. The resulting loop based at $x_{i}$ is null-homotopic in $X$,
and hence so is $\gamma $.
\end{proof}

\begin{lemma}
\label{315conn}Let $X$ be a chain connected uniform space. If there is an
entourage $E$ of $X$ such that the $E$-balls of $X$ are connected (resp.
pathwise connected) then $X$ is connected (resp. pathwise connected).
\end{lemma}

\begin{proof}
It will follow from\ Proposition \ref{consume} and a standard theorem from
topology if we show by induction that all $E^{n}$-balls are connected. The $%
n=1$ case is given. Suppose that $E^{n}$-balls and $E$-balls are all
conneced for some $n$. Let $y\in B(x,E^{n+1})$. By definition of $%
B(x,E^{n+1})$ there is some $z$ such that $z\in B(x,E^{n})$ and $z\in B(y,E)$%
. Each of these balls is connected and they intersect in $z$, hence their
union is connected set. We have shown that every element of $B(x,E^{n+1})$
is contained in a connected subset of $B(x,E^{n+1})$ containing $x$, and so $%
B(x,E^{n+1})$ is connected. If the $E$-balls are all pathwise connected then
by what we have just proved $X$ is connected and locally pathwise connected,
hence pathwise connected.
\end{proof}

\begin{corollary}
If $X$ is chain connected and uniformly locally connected (resp. uniformly
locally pathwise connected) then $X$ is connected (resp. pathwise connected).
\end{corollary}

\begin{proposition}
\label{52}Let $X$ be a uniform space and $E$ be an entourage in $X$ such
that the $E$-balls have one of the following properties: chain connected,
connected, pathwise connected. Then $X_{E}$ has the same property.
\end{proposition}

\begin{proof}
Suppose that the $E$-balls of $X$ are chain connected. Let $F\subset E$ be
an entourage and $\alpha =\{\ast =x_{0},...,x_{n}\}$ be an $E$-chain. Since
the $E$-balls are chain connected we may suppose up to $E$-homotopy that $%
\alpha $ is an $F$-chain. For example, we may join $x_{0}$ and $x_{1}$ by an 
$F$-chain $\{x_{0},y_{1},...,y_{m},x_{1}\}$ that lies entirely in $%
B(x_{0},E) $. We may remove the points $y_{1},..y_{m}$ one at a time in
order to obtain an $E$-homotopy from $%
\{x_{0},y_{1},...,y_{m},x_{1},...,x_{n}\}$ to $\alpha $. Then letting $%
\alpha _{i}:=\{\ast =x_{0},...,x_{i}\}$ we have that $([\alpha
_{i}]_{E},[\alpha _{i+1}]_{E})\in F^{\ast }$ for all $i$. In particular, $%
\{[\alpha _{0}]_{E},...,[\alpha _{n}]_{E}=[\alpha ]_{E}\}$ is an $F^{\ast }$%
-chain to $[\ast ]_{E}$. This shows that $X_{E}$ is chain connected. Now if $%
E$ has connected (resp. pathwise connected) balls then these balls are chain
connected and therefore $X_{E}$ is chain connected by what we just proved.
According to Proposition \ref{xee}, the $E^{\ast }$-balls in $X_{E}$ are
connected (resp. pathwise connected), and Lemma \ref{315conn} now shows that 
$X_{E}$ is connected (resp. pathwise connected).
\end{proof}

\begin{theorem}
\label{0627simply}If $X$ is a uniformly locally pathwise connected,
connected and simply connected uniform space then $X$ is universal.
\end{theorem}

\begin{proof}
Let $E$ be an entourage with pathwise connected open balls. Then $X_{E}$ is
pathwise connected by Proposition \ref{52}. Moreover, by Corollary 34 of 
\cite{P}, the action of $\delta _{E}(X)$ on $X_{E}$ is properly
discontinuous. Since $X$ is a Poincar\'{e} space, $\phi
_{XE}:X_{E}\rightarrow X$ is a traditional cover (cf. \cite{M}) and since $X$
is simply connected, $\phi _{XE}$ must be a trivial cover, hence bijective.
This means that the bi-uniformly continuous mapping $\phi _{XE}$ is a
uniform homeomorphism.
\end{proof}

\begin{corollary}
\label{compactu}Every compact, connected, locally pathwise connected, simply
connected topological space is universal.
\end{corollary}

Proposition \ref{liftlem} now implies:

\begin{corollary}
\label{pathlift}If $c:[0,1]\rightarrow X$ is a (continuous) path, where $X$
is uniform, with $c(0)=\ast $, then there is a unique lift $%
c_{L}:[0,1]\rightarrow \widetilde{X}$ such that $c_{L}(0)=\ast $ and $\phi
\circ c_{L}=c$, where $\phi :\widetilde{X}\rightarrow X$ is the projection.
A similar statement holds for homotopies.
\end{corollary}

If $c:[0,1]\rightarrow X$ is a path from $\ast $ to $x\in X$ then $x=\phi
\circ c_{L}(1)\in \phi (\widetilde{X})$. We have:

\begin{corollary}
\label{pathsur}If $X$ is a uniform space with projection $\phi :\widetilde{X}%
\rightarrow X$ then the pathwise connected component of $X$ at $\ast $ is
contained in $\phi (P)$, where $P$ is the pathwise connected component of $%
\ast $ in $\widetilde{X}$.
\end{corollary}

\begin{corollary}
\label{0628new}If $X$ is a connected, uniformly locally pathwise connected
uniform space then $X$ is coverable. In particular, every compact,
connected, locally pathwise connected topological space is coverable.
\end{corollary}

\begin{proof}
We may index the fundamental system of $X$ using entourages having pathwise
connected balls. But for each such $E$, $X_{E}$ is pathwise connected by
Proposition \ref{52} and therefore each $\phi _{E}:\widetilde{X}\rightarrow
X_{E}$ is surjective by the previous corollary.
\end{proof}

\begin{theorem}
\label{0518simp}Suppose $X$ is a uniformly locally pathwise connected
uniform space having an entourage $E$ such that for all $x\in X$ every loop
in $B(x,E)$ is null-homotopic in $X$ and an entourage $F$ having pathwise
connected balls such that $F^{2}\subset E$. Then $X_{F}$ is pathwise
connected, uniformly locally pathwise connected and simply connected, hence
universal.
\end{theorem}

\begin{proof}
According to Proposition \ref{52}, $X_{F}$ is pathwise connected, and the $%
F^{\ast }$-balls, which are uniformly homeomorphic to the $F$-balls, are
pathwise connected. By Proposition \ref{identify} the function $\phi
_{F^{\ast }X_{F}}:\left( X_{F}\right) _{F^{\ast }}\rightarrow X_{F}$ is a
uniform homeomorphism. This implies that every $F^{\ast }$-loop in $X_{F}$
is $F^{\ast }$-homotopic to the trivial loop. Now let $c:[0,1]\rightarrow
X_{F}$ be a loop based at $\ast $. Since $c$ is uniformly continuous there
is some $F^{\ast }$-loop $\gamma =\{\ast =x_{0},...,x_{n}=\ast \}$ such that
each $x_{i}=c(t_{i})$ for some $i$ and $c([t_{i},t_{i+1}])\subset
B(x_{i},F^{\ast })$ for all $i$. We may use any $F^{\ast }$-homotopy from $%
\gamma $ to $\ast $ to construct a null-homotopy from $c$ to the trivial
loop as follows. Suppose that $(x_{i},y),(x_{i+1},y)\in F^{\ast }$ for some $%
y\in X$. We may join $x_{i}$ to $y$ and $y$ to $x_{i+1}$ by paths $\alpha
_{1}$ and $\alpha _{2}$ in $B(y,F^{\ast })$. Let $\alpha _{3}$ be the
restriction of $c$ to $[t_{i},t_{i+1}]$. Then we have a loop $\beta =\alpha
_{1}\ast \alpha _{2}\ast \alpha _{3}^{-1}$, each segment of which lies in
the $F^{\ast }$-ball centered at its endpoint. Now let $\psi :=\phi
_{XF}\circ \beta $, which is a loop consisting of three segments, each of
which lies in the $F$-ball centered at its endpoint. Since $F^{2}\subset E$, 
$\psi $ lies entirely in the $E$-ball centered at a point, and hence is
null-homotopic in $X$. According to Theorem \ref{76} any null-homotopy of $%
\psi $ lifts to a null-homotopy of $\beta $. In other words, the path $c$ is
homotopic to a path $c_{1}$ obtained by replacing $\alpha _{3}$ by $\alpha
_{1}$ concatenated by $\alpha _{2}$. We may carry out a similar process when
removing any point in $\gamma $ to form a new $F^{\ast }$-chain as part of
an $F^{\ast }$-homotopy of $\gamma $. Finitely many such steps show that $c$
is null-homotopic. Since the $F$-balls of $X$, and hence the $F^{\ast }$%
-balls of $X_{F}$, are pathwise connected, $X_{F}$ is uniformly locally
pathwise connected and pathwise connected, hence universal by Theorem \ref%
{0627simply}.
\end{proof}

\section{The homomorphism $\protect\lambda :\protect\pi _{1}(X)\rightarrow 
\protect\delta _{1}(X)$}

Let $c:[0,1]\rightarrow X$ be a path in a coverable space $X$ such that $%
c(0)=c(1)=\ast $. Since $c(1)=\ast $, $\phi (c_{L}(1))=\ast $, where $\phi $
is the universal covering map of $X$ and $c_{L}$ is the lift given by
Proposition \ref{liftlem}. That is, $c_{L}(1)\in \delta _{1}(X)$. Define a
function $\lambda :\pi _{1}(X)\rightarrow \delta _{1}(X)$ by $\lambda
([c])=c_{L}(1)$, where $[c]$ is the homotopy equivalence class of a loop $c$
based at $\ast $. This function is well defined. In fact, suppose $c$ and $d$
are loops based at $\ast $ that are homotopic via a homotopy $\eta
:[0,1]\times \lbrack 0,1]\rightarrow X$, with $\eta (0,t)=c(t)$ and $\eta
(1,t)=d(t)$. Then $\eta _{L}(t,1)$ is a path joining $c_{L}(1)$ to $d_{L}(1)$%
, both of which lie in $\phi ^{-1}(\ast )$. But $\phi ^{-1}(\ast )=\delta
_{1}(X)$ is a prodiscrete group (hence totally disconnected) with respect to
the subspace topology. Therefore $c_{L}(1)=d_{L}(1)$.

\begin{proposition}
\label{lifter}Let $X$ be coverable, $[c]\in \pi _{1}(X)$, where $%
c:[0,1]\rightarrow X$ is based at $\ast $, and $E$ be an entourage in $X$.
If $0=t_{0}<\cdot \cdot \cdot <t_{n}=1$ are such that for every $i$, $%
c(t)\in B(c(t_{i}),E))$ for all $t\in \lbrack t_{i},t_{i+1}]$ then $\theta
_{E}(\lambda ([c]))=[\ast =c(t_{0}),...,c(t_{n})]_{E}$. (Such values $t_{i}$
always exist since $c$ is uniformly continuous.)
\end{proposition}

\begin{proof}
Note that by definition,%
\begin{equation*}
\theta _{E}(\lambda ([c]))=\theta _{E}(c_{L}(1))=c_{E}(1)
\end{equation*}%
where $c_{E}$ is the lift given by Proposition \ref{liftlem}. Let $\gamma
=\{\ast =c(t_{0}),...,c(t_{n})\}$, where the points $t_{i}$ satisfy the
conditions of the proposition. We will show by induction on $i$ that $%
c_{E}(t_{i})=[c(t_{0}),c(t_{1}),...,c(t_{i})]_{E}$ for all $0\leq i\leq n$.
The case $i=0$ is trivial. Suppose that we have proved the statement for
some $i$, $0\leq i<n$. By definition, 
\begin{equation*}
([c(t_{0}),c(t_{1}),...,c(t_{i})]_{E},[c(t_{0}),c(t_{1}),...,c(t_{i}),c(t_{i+1})]_{E})\in E^{\ast }
\end{equation*}%
and $\phi _{XE}([c(t_{0}),c(t_{1}),...,c(t_{i}),c(t_{i+1})]_{E})=c(t_{i+1})$%
. Recall from Proposition \ref{xee} that $\phi _{XE}$ is a uniform
homeomorphism when restricted to any $E^{\ast }$-ball of $X_{E}$, and since
we also have $\phi _{XE}(c_{E}(t_{i+1}))=c(t_{i+1})$ we will be finished if
we can show 
\begin{equation*}
([c(t_{0}),c(t_{1}),...,c(t_{i})]_{E},c_{E}(t_{i+1}))\in E^{\ast }\text{.}
\end{equation*}%
But the unique lift $c_{i}$ of $c\mid _{\lbrack t_{i},t_{i+1}]}$ to $X_{E}$
starting at $c_{E}(t_{i})$ must be $\eta \circ c\mid _{\lbrack
t_{i},t_{i+1}]}$, where $\eta $ is the inverse of $\phi _{XE}$ restricted to 
$B([c(t_{0}),c(t_{1}),...,c(t_{i})]_{E},E^{\ast })$. Since $c_{E}\mid
_{\lbrack t_{i},t_{i+1}]}$ also satisfies the conditions for this lift, we
must have that 
\begin{equation*}
c_{E}(t_{i+1})=c_{i}(t_{i+1})=\eta \circ c(t_{i+1})\in
B([c(t_{0}),c(t_{1}),...,c(t_{i})]_{E},E^{\ast })\text{.}
\end{equation*}
\end{proof}

\begin{theorem}
\label{nat}If $X$ is coverable then the natural mapping $\lambda :\pi
_{1}(X)\rightarrow \delta _{1}(X)$ is a homomorphism, the image of which is
the normal subgroup $\sigma (X)$ of $\delta _{1}(X)$ that leaves invariant
the pathwise connected component of $\widetilde{X}$ containing $\ast $. In
particular, $\lambda $ is surjective if and only if $\widetilde{X}$ is
pathwise connected.
\end{theorem}

\begin{proof}
To see why $\lambda $ is a homomorphism note that $\phi \circ (c_{L}\ast
\left( \lambda ([c])\circ d_{L}\right) )=c\ast d$ (here $\lambda ([c])\circ
d_{L}$ is the \textquotedblleft translate\textquotedblright\ via the uniform
homeomorphism $\lambda ([c])$ of $d_{L}$ to the endpoint of $c_{L}$ and $%
\phi $ is the universal covering map of $X$) and therefore by uniqueness, $%
c_{L}\ast \left( \lambda ([c])\circ d_{L}\right) =\left( c\ast d\right) _{L}$%
. But the endpoint of $c_{L}\ast \lambda ([c])\circ d_{L}$ is $%
c_{L}(1)d_{L}(1)=\lambda ([c])\lambda ([d])$. If $g\in \delta _{1}(X)$
stabilizes the pathwise connected component of $\widetilde{X}$ containing $%
\ast $ then $\ast $ and $g(\ast )$ are joined by an path $\alpha $, and $%
\phi (\alpha )$ is a loop based at $\ast $ with $\lambda ([\alpha ])=g$. On
the other hand, if $g$ lies in the image of $\lambda $ then $\ast $ and $%
g(\ast )$ are joined by an path (namely the lift of a loop, the image of
whose equivalence class is $g$). Now suppose $x$ lies in the pathwise
connected component of $\widetilde{X}$ containing $\ast $. We may join $x$
to $\ast $ by an path $c$. But then $g\circ c$ joins $g(x)$ and $g(\ast )$,
and we have already observed that the latter is joined to $\ast $ by an path.
\end{proof}

The kernel of $\lambda $, by definition, consists of those $[c]$ such that $%
c $ lifts to a loop in $\widetilde{X}$, or equivalently the lift of $c$ to
each $X_{E}$ is a loop. It is of obvious interest when $\ker \lambda $ is
trivial. This can be checked in some special cases. For example, in the case
of the Hawaiian earring, a problem that essentially amounts to the
injectivity of $\lambda $ occupied several papers (\cite{G}, \cite{MM}, \cite%
{D}) and was simply stated in \cite{B} (along with an incorrect statement
that amounts to surjectivity of $\lambda )$.

\begin{proposition}
\label{inj}Let $X$ be coverable. Then

\begin{enumerate}
\item If $\widetilde{X}$ is simply connected then $\lambda :\pi
_{1}(X)\rightarrow \delta _{1}(X)$ is injective.

\item If $X$ is pathwise connected and $\lambda :\pi _{1}(X)\rightarrow
\delta _{1}(X)$ is injective then $\widetilde{X}$ is simply connected.
\end{enumerate}
\end{proposition}

\begin{proof}
If $\widetilde{X}$ is simply connected and $[c]\in \ker \lambda $ then $c$
lifts to a loop $c_{L}$ in $\widetilde{X}$ based at $\ast $, which is then
null-homotopic. The image of any null-homotopy of $c_{L}$ via the universal
covering map $\phi $ is a null-homotopy of $c$.

Now suppose that $X$ is pathwise connected and $\lambda :\pi
_{1}(X)\rightarrow \delta _{1}(X)$ is injective. Let $c:[0,1]\rightarrow 
\widetilde{X}$ be a loop based at some point $x$. Then $\phi \circ c$ is a
loop in $X$, where $\phi $ is the universal covering map of $X$. Let $%
d:[-1,2]\rightarrow X$ be the concatenation of a path $\alpha $ from $\ast $
to $\phi (c(0))$, followed by $\alpha ^{-1}$, parameterized so that the
restriction to $[0,1]$ is $\phi \circ c$. So $d$ represents an element of $%
\pi _{1}(X)$. Let $g\in \delta _{1}(X)$ be such that $g(d_{L}(0))=c(0)$
(such a $g$ exists since $\phi (c(0))=d(0)$). Now the composition of $g$
with $d_{L}\mid _{\lbrack 0,1]}$ is simply $c$ and therefore $d_{L}\mid
_{\lbrack 0,1]}$ is a loop. But then $\lambda ([d])=0$ and since $\lambda $
is injective, $d$ must be null-homotopic. But then any null-homotopy of $d$
lifts to one of $d_{L}$, and the composition of the lifted homotopy with $g$
gives rise to a null-homotopy of $c$.
\end{proof}

\begin{proposition}
\label{dense}If $X$ is a uniformly locally pathwise connected, connected
uniform space then the pathwise connected component of $\widetilde{X}$ is
dense in $\widetilde{X}$ and $\lambda (\pi _{1}(X))$ is dense in $\delta
_{1}(X)$.
\end{proposition}

\begin{proof}
These two statements follow from the following general result concerning
induced functions on inverse limits: Let $(X_{\alpha },\phi _{\alpha \beta
}) $ be an inverse system of topological spaces with continuous bonding maps
and $X:=\underleftarrow{\lim }X_{\alpha }$. Let $f_{\alpha }:Y\rightarrow
X_{\alpha }$ be a collection of continuous surjections from a topological
space $Y$ such that $\phi _{\alpha \beta }\circ f_{\beta }=f_{\alpha }$ for
all $\alpha \leq \beta $. Then the induced mapping $f:Y\rightarrow X$ has
dense image in $X$. We do not have a reference for this exact statement but
the proof is straightforward and similar to the proof of III.7.3 Proposition
2 in \cite{Bk}. Now if $E$ has pathwise connected balls then $X_{E}$ is
pathwise connected and by Corollary \ref{pathsur} the restriction $\psi $ of 
$\phi _{E}$ to the pathwise connected component $P$ of $\widetilde{X}$ is
surjective onto $X_{E}$ and the proof of the first part is finished by the
above general statement.

On the other hand, if $x\in \delta _{E}(X)$ then there is some path $\alpha $
from $\ast $ to $x$, and $\phi _{XE}\circ \alpha $ is a loop $\gamma $ in $X$
based at $\ast $ such that the unique lift $\gamma _{E}$ has $x$ as its
endpoint. But then $\phi _{E}\circ \lambda ([\gamma ])=x$. In other words, $%
\theta _{E}\circ \lambda $ is surjective, and the proof is finished by the
above general statement.
\end{proof}

Note that the above proof really only requires that there be a basis for the
uniformity of $X$ such that for each $E$ in the basis, $X_{E}$ is pathwise
connected. Since the closure of a connected set is connected we have:

\begin{corollary}
If $X$ is a uniformly locally pathwise connected, connected uniform space
then $\widetilde{X}$ is connected.
\end{corollary}

\begin{proposition}
\label{newquotient}If $X$ is a locally pathwise connected and connected
space and $P$ is the pathwise connected component of $\widetilde{X}$
containing $\ast $ then $X$ is the quotient of $P$ by the free isomorphic
action of $\lambda (\pi _{1}(X)$.
\end{proposition}

\begin{proof}
We already know that $\lambda (\pi _{1}(X))$ acts freely and isomorphically
on $P$ since $\delta _{1}(X)$ does. We need to check that the restriction $%
\psi $ of $\phi $ to $P$ is bi-uniformly continuous and that for any $x\in P$%
, $\psi ^{-1}(\psi (x))$ is precisely the orbit $\lambda (\pi _{1}(X))(x)$
of $x$ (see Remark 13 of \cite{P}). We know that $\psi $ is uniformly
continuous. Let $D:=\phi _{E}^{-1}(F^{\ast })\cap P$ be an entourage in $P$;
we may assume that the $E$-balls are pathwise connected, hence $X_{E}$ is
pathwise connected, hence the restriction $\psi _{E}$ of $\phi _{E}$ to $P$
is surjective by Corollary \ref{pathsur}. Since $\psi _{E}$ is surjective, 
\begin{equation*}
F^{\ast }=\psi _{E}(\psi _{E}^{-1}(F^{\ast }))=\psi _{E}(\phi
_{E}^{-1}(F^{\ast })\cap P)=\psi _{E}(D)\text{.}
\end{equation*}%
Now 
\begin{equation*}
\psi (D)=\phi (D\cap P)=\phi _{XE}(\phi _{E}(D\cap P))=\phi _{XE}(\psi
_{E}(D))=\phi _{XE}(F^{\ast })
\end{equation*}%
which is an entourage since $\phi _{XE}$ is bi-uniformly continuous. This
shows that $\psi $ is bi-uniformly continuous. The statement about the
orbits simply follows from the fact that the orbits of $\psi $ are precisely
the orbits of $\phi $ intersected with $P$ and that $\psi ^{-1}(y)=\phi
^{-1}(y)\cap P$ for any $y\in X$.
\end{proof}

\begin{remark}
Note that when $\lambda $ is injective, the action in the above proposition
is in fact an action by $\pi _{1}(X)$. On the other hand, when $\lambda $ is
not surjective one has the disadvantage that $P$ is not complete with the
uniformity induced by $\widetilde{X}$.
\end{remark}

\begin{definition}
A uniform space $X$ is called strongly coverable if $X$ is chain connected
and for some entourage $E$, $X_{E}$ is universal, hence the universal cover
of $X$.
\end{definition}

It is clear from the discussion in the introduction that the Topologist's
Sine Curve is strongly coverable.

\begin{lemma}
\label{strong}If $X$ is strongly coverable and uniformly locally pathwise
connected then $\lambda :\pi _{1}(X)\rightarrow \delta _{1}(X)$ is
surjective.
\end{lemma}

\begin{proof}
Let $E$ be an entourage such that $X_{E}$ is universal and $F\subset E$ be
an entourage having pathwise connected balls. Then $X_{F}$ is pathwise
connected and $F^{\ast }$ is in the universal base of $X_{E}$ by Corollary %
\ref{0518u}. Hence $X_{E}$ and $X_{F}$ are uniformly homeomorphic and
therefore $X_{E}=\widetilde{X}$ is also pathwise connected. The proof is
finished by Theorem \ref{nat}.
\end{proof}

From Theorem \ref{0518simp}, Lemma \ref{strong}, and Proposition \ref{inj}
we obtain:

\begin{theorem}
Every connected, uniformly locally pathwise connected and uniformly
semi-locally simply connected uniform space $X$ is strongly coverable with
pathwise connected, simply connected universal cover and $\delta _{1}(X)=\pi
_{1}(X)$.
\end{theorem}

\begin{corollary}
\label{0518equiv}If $X$ is a compact Poincar\'{e} space then $\widetilde{X}$
is the traditional universal cover of $X$ and $\lambda :\pi
_{1}(X)\rightarrow \delta _{1}(X)$ is an isomorphism.
\end{corollary}

\section{Dimension and the universal cover}

We refer the reader to \cite{I} for more background on dimension and uniform
spaces. Suppose that $X$ has uniform dimension $\leq n$, which we denote by $%
u\dim X\leq n$. This means that any uniform open cover $\mathcal{V}$ of $X$
has a refinement by a uniform open cover of order $n+1$ (a uniform open
cover is an open cover that is refined by the cover of $X$ by $F$-balls for
some fixed entourage $F$). Note that this particular notion of dimension is
called \textquotedblleft large dimension\textquotedblright\ in \cite{I}, and
is denoted by $\triangle dX$. We will use theorems from \cite{I} concerning
another dimension, called \textquotedblleft uniform
dimension\textquotedblright\ in \cite{I} and denoted by $\delta dX$. This
particular notion of dimension uses finite covers and is somewhat more
difficult to work with in the present situation. However, it is always true
that $\delta dX\leq \triangle dX=u\dim X$ (Theorem V.5 in \cite{I}--in fact
the dimensions are equal if $\triangle dX$ is finite) and therefore if we
know that $u\dim X\leq n$ then we may use theorems from \cite{I} that
require $\delta dX\leq n$. If $X$ is compact then both of these dimensions
are equal to covering dimension, which we denote by $\dim X$.

\begin{proposition}
\label{dimprop}If $X$ is a uniform space with $u\dim X\leq n$ then for any
entourage $E$, $u\dim X_{E}\leq n$.
\end{proposition}

\begin{proof}
Every uniform open cover of $X_{E}$ is refined by the cover of $X_{E}$ by $%
F^{\ast }$-balls for some entourage $F$ such that $F^{2}\subset E$ and
therefore we need only consider the open cover of $X_{E}$ by $F^{\ast }$%
-balls for such $F$. By definition of uniform dimension, the cover of $X$ by 
$F$-balls has a refinement by a uniform open cover $\mathcal{V}$ such that
every $x\in X$ is contained in at most $n+1$ sets in $\mathcal{V}$. Let $%
A\in \mathcal{V}$. Then $A\subset B(x,F)$ for some $x\in X$. Let $\mathcal{W}%
_{A}$ be the collection of all sets of the form $\phi _{XE}^{-1}(A)\cap
B(y,F^{\ast })$ where $y\in \phi _{XE}^{-1}(x)$ and $\mathcal{W}%
:=\bigcup\limits_{A\in \mathcal{V}}\mathcal{W}_{A}$. First, $\mathcal{W}$ is
a cover. In fact, if $z\in X_{E}$, $\phi _{XE}(z):=w\in A$ for some $A\in 
\mathcal{V}$, with $A\subset B(x,F)$. But $\phi _{XE}$ restricted to $%
B(z,E^{\ast })$ is a uniform homeomorphism onto $B(w,E)$, which contains $%
A\subset B(x,F)$ since $F^{2}\subset E$. Therefore there is some $y\in
B(z,E^{\ast })\cap \phi _{XE}^{-1}(x)$. But $\phi _{XE}(B(z,F^{\ast
}))=B(w,F)$ and $(w,x)\in F$; therefore $(z,y)\in F^{\ast }$. That is, $z\in
\phi _{XE}^{-1}(A)\cap B(y,F^{\ast })$. Since $\mathcal{V}$ is a uniform
cover, so is $\mathcal{W}$. In fact, if the $D$-ball cover of $X$ refines $%
\mathcal{V}$ for some $D\subset F$ then the $D^{\ast }$-ball cover refines $%
\mathcal{W}$. By definition, $\mathcal{W}$ refines the $F^{\ast }$-ball
cover of $X_{E}$. Finally, to check that the order of $\mathcal{W}$ is at
most $n+1$, we need only check that if $y,z\in \phi _{XE}^{-1}(x)$ are
distinct and $A\in \mathcal{V}$ then 
\begin{equation*}
\left[ \phi _{XE}^{-1}(A)\cap B(y,F^{\ast })\right] \cap \left[ \phi
_{XE}^{-1}(A)\cap B(z,F^{\ast })\right] =\varnothing \text{.}
\end{equation*}%
But $B(y,F^{\ast })$ and $B(z,F^{\ast })$ are already disjoint, because
otherwise for any $w\in B(y,F^{\ast })\cap B(z,F^{\ast })$ we would have $%
z\in B(y,\left( F^{\ast }\right) ^{2})\subset B(y,E^{\ast })$, which
contradicts the fact that $\phi _{XE}$ is injective on $E^{\ast }$-balls.
\end{proof}

In \cite{I}, Theorem IV.32 it is shown that the inverse limit of spaces with
uniform dimension at most $n$ must have uniform dimension at most $n$.
Therefore by Proposition \ref{dimprop}:

\begin{theorem}
\label{dimthm}If $X$ is a coverable uniform space with $u\dim X\leq n$ then $%
u\dim \widetilde{X}\leq n$.
\end{theorem}

Note that essentially the same argument shows that if $X$ and $Y$ are
uniform spaces, $f:X\rightarrow Y$ is a cover, and $u\dim Y\leq n$ then $%
u\dim X\leq n$. We conjecture that in this situation $X$ and $Y$ have
exactly the same dimension.

\begin{proposition}
\label{circle} If $X$ is coverable uniform space with $u\dim X\leq 1,$ then $%
\tilde{X}$ contains no simple closed curve (i.e. topological circle).
\end{proposition}

\begin{proof}
Suppose that there is a topological embedding $f:S^{1}\rightarrow S\subset 
\tilde{X}$ with the inverse homeomorphism $g:S\rightarrow S^{1}$. By Theorem %
\ref{dimthm} we have $u\dim \widetilde{X}\leq 1$ and therefore by Theorem
V.13 of \cite{I} there is an extension of $g$ to a uniformly continuous
function $G:\tilde{X}\rightarrow S^{1}$. Choose any point $\ast $ as the
basepoint in both $S$ and $\widetilde{X}$, and choose $g(\ast )$ as the
basepoint in $S^{1}$ (see Remark \ref{base2}). Let $\psi :\mathbb{R}=\tilde{%
S^{1}}\rightarrow S^{1}$ be the traditional universal cover of $S^{1}$ also
with some choice of basepoint $\ast $ in $\psi ^{-1}(\ast ),$ which is also
the universal cover in the sense of the present paper by Corollary \ref%
{0518equiv}. Since $\tilde{X}$ is universal by Theorem \ref{xic},
Proposition \ref{liftlem} (2) implies that there is a unique lift $G_{L}:%
\tilde{X}\rightarrow \mathbb{R}$ such that $G_{L}(\ast )=\ast $ and $\psi
\circ G_{L}=G.$ Thus we get that 
\begin{equation*}
\psi \circ G_{L}\circ f=G\circ f=g\circ f=id_{S^{1}}\text{.}
\end{equation*}%
This implies that $G_{L}\circ f$ is a topological embedding of $S^{1}$ into $%
\mathbb{R}$, which is impossible.
\end{proof}

\begin{theorem}
\label{laminj}If $X$ is coverable uniform space with $u\dim X\leq 1$ (in
particular if $X$ is compact with covering dimension $\dim X\leq 1$), then $%
\widetilde{X}$ is simply connected and if $X$ is pathwise connected the
homomorphism $\lambda :\pi _{1}(X,\ast )\rightarrow \delta _{1}(X)$ is
injective.
\end{theorem}

\begin{proof}
By Proposition \ref{inj}, it is enough to prove that $\tilde{X}$ is simply
connected. Suppose that $c:[0,1]\rightarrow \tilde{X}$ is a loop in $\tilde{X%
}$. Then its image $C$ is a Peano continuum that contains no simple closed
curves by Proposition \ref{circle}. By the Hahn-Muzurkiewicz Theorem, $C$ is
locally connected, hence a dendrite (see section 51, VI in \cite{Ku}). Then $%
C$ is contractible by Corollary 7 in Section 54, VII of \cite{Ku}. Therefore 
$C$ has trivial fundamental group. This means that any loop at any basepoint
in $C$ (including the loop $c$) is null-homotopic in $C$ hence in $%
\widetilde{X}$.
\end{proof}

Combining the above theorem with Proposition \ref{newquotient} we have the
following:

\begin{corollary}
Let $X$ be a pathwise connected, uniformly locally pathwise connected
uniform space with $u\dim X=1$. Then $X$ is the quotient of a
one-dimensional pathwise connected, simply connected uniform space via a
free isomorphic action of $\pi _{1}(X)$.
\end{corollary}

\section{Pseudometric spaces\label{pssec}}

\begin{definition}
\label{TCC}Let $X$ be a uniform space. We define an entourage $E$ to be
chain connected if every $E$-ball in $X$ is chain connected. We say $X$ is
totally chain connected if $X$ has a uniformity base that includes $X\times
X $ such that each entourage in the base is chain connected.
\end{definition}

From Proposition \ref{52} we immediately have:

\begin{lemma}
\label{315totally}If $X$ is a totally chain connected uniform space and $E$
is a chain connected entourage then $X_{E}$ is totally chain connected. If
moreover $X$ is uniformly locally connected (resp. uniformly locally
pathwise connected) then $X_{E}$ has the same property.
\end{lemma}

\begin{example}
In the proof of Theorem 9, \cite{BPLCG}, it was shown that the character
group of $\mathbb{Z}^{\mathbf{N}}$ is not coverable. But this group is known
to be connected and locally connected, hence totally chain connected.
\end{example}

\begin{theorem}
Every totally chain connected pseudometric space $X$ is coverable.
\end{theorem}

\begin{proof}
Since $X$ is a pseudometric space we can find a countable sequence $%
\{E_{i}\} $ of chain connected entourages forming a base for the uniformity
of $X$ (hence cofinal in the set of all entourages of $X$). Therefore $%
\widetilde{X}=\underleftarrow{\lim }X_{E_{i}}$. Each of the spaces $%
X_{E_{i}} $ is chain connected by Lemma \ref{315totally} and when $j\geq i$
we may identify $X_{E_{j}}$ with $\left( X_{E_{i}}\right) _{E_{j}^{\ast }}$
by Proposition \ref{identify}. Therefore, Lemma \ref{lemsur5} implies that
each bonding map $\phi _{E_{i}E_{j}}$ is surjective. Since the inverse
system is countable it follows that the projections $\phi _{E_{i}}:%
\widetilde{X}\rightarrow X_{E_{i}}$ are all surjective.
\end{proof}

Since connected sets are chain connected, we have the following two
corollaries:

\begin{corollary}
\label{pseu}Every connected, uniformly locally connected pseudometric space
is coverable.
\end{corollary}

Recall that a Peano continuum is a Hausdorff topological space that is the
continuous image of an interval. Equivalently (by the Hahn-Mazurkiewicz
Theorem), a Peano continuum is a compact, connected, locally (pathwise)
connected metrizable space.

\begin{corollary}
\label{peano}Every Peano continuum $X$ is coverable. Moreover, there is a
compact subset $S$ of the pathwise connected component of $\widetilde{X}$
such that the restriction of the universal covering map $\phi $ to $S$ is
surjective.
\end{corollary}

\begin{proof}
The first statement follows from Corollary \ref{0628new}. Let $%
c:[0,1]\rightarrow X$ be a continuous surjection. Then $S:=c_{L}([0,1])$ has
the desired properties.
\end{proof}

\begin{notation}
\label{metnote}To simplify matters, when $X$ is a metric space we refer to
an $E_{\varepsilon }$-loop, where $E_{\varepsilon }$ is the metric entourage
having open $\varepsilon $-balls as $E_{\varepsilon }$-balls, as an $%
\varepsilon $-loop ($\varepsilon >0$), $E_{\varepsilon }$-homotopies as $%
\varepsilon $-homotopies, etc. We will denote $X_{E_{\varepsilon }}$ by $%
X_{\varepsilon }$, $\phi _{E_{\varepsilon }E_{\delta }}$ by $\phi
_{\varepsilon \delta }$, $\phi _{XE_{\varepsilon }}$ by $\phi _{X\varepsilon
}$, and $\phi _{E_{\varepsilon }}$ by $\phi _{\varepsilon }$. 
\end{notation}

Note that for a metric space $X$ and $\varepsilon >0$, the cover $\phi
_{X\varepsilon }:X_{\varepsilon }\rightarrow X$ is a broadening of the
notion of \textquotedblleft $\varepsilon $-cover\textquotedblright ,
something that goes back at least to Spanier's book (\cite{S}), and has been
used in \cite{SW1} and \cite{SW2} to study universal covers of limits of
Riemannian manifolds. However, the construction of $\varepsilon $-covers
uses paths and standard homotopies rather than chains.

It is not hard to construct examples of uniformly locally pathwise connected
metric spaces having metric balls that are not necessarily pathwise
connected, and so that the metric entourages are not covering entourages.
Recall that an \textit{inner metric space} is a metric space such that the
distance between any two points is the infimum of lengths of curves joining
them. A \textit{geodesic} space further has the property that the distance
is realized as the length of some curve, called a minimal geodesic. The
metric balls of any inner metric space are pathwise connected and therefore
the metric entourages are covering entourages. We have:

\begin{corollary}
\label{inner}Every inner metric space $X$ is coverable. Moreover, $%
\widetilde{X}=\underleftarrow{\lim }X_{\varepsilon _{i}}$, where $%
\varepsilon _{i}$ is any sequence of positive values decreasing to $0$.
\end{corollary}

From Lemma \ref{315totally} and Corollary \ref{uit} we have:

\begin{corollary}
An inner metric space $X$ is universal if and only if for any $\varepsilon
>0 $, every $\varepsilon $-loop based at some point is $\varepsilon $%
-homotopic to the trivial loop.
\end{corollary}

Note that the topologist's sine curve, discussed in the introduction, is
coverable but not totally chain connected because at some points every small
neighborhood contains pairs of points that cannot be joined by an
arbitrarily fine chain that stays inside the neighborhood.

We next give a unified calculation of the deck groups $\delta _{1}(X)$ of
the Hawaiian earring $X=H$, the Sierpin'ski gasket $X=\Delta _{S}$, and
Sierpin'ski carpet $X=C_{S}$.

This unified calculation is possible because each of these spaces admits a
description as a countable intersection 
\begin{equation}
X=\bigcap_{n=0}^{\infty }X_{n}.  \label{inters}
\end{equation}

Let us begin with the case $X=H.$ Define $X_{0}:=H_{0}$ as the closed disc $%
D $ of diameter 1, lying in the upper Euclidean half-plane $y\geq 0$ and
tangent to the $x$-axis at the origin $O=(0,0)$. Let $C=\partial D$ be its
boundary circle. Denote by $\psi _{\{\alpha ,c\}}$ the $\alpha $-homothety
of Euclidean plane $E^{2}$ with the center $c.$ We omit $c$ in the notation
if $c=O.$ Now define $X_{n}$, $n\geq 0,$ iteratively by the formulas 
\begin{equation}
X_{0}:=D\text{,}\quad X_{n+1}:=\psi _{\frac{1}{2}}(X_{n})\cup C\text{.}
\label{ind1}
\end{equation}%
Then $X=H$ is defined by Formulas (\ref{ind1}) and (\ref{inters}).

For $X=\Delta _{S}$ let $X_{0}:=\Delta _{0}$ be the isosceles triangle of
diameter 1 in the Euclidean quarter-plane $x\geq 0,y\geq 0$, with one side
on the $x$-axis and one vertex $v_{0}=O$; denote the other vertices by $%
v_{1},v_{2}$. Now define $X_{n}$ iteratively by the formulas 
\begin{equation}
X_{0}=\Delta _{0},\quad X_{n+1}:=\bigcup_{i=0}^{2}\psi _{\{\frac{1}{2}%
,v_{i}\}}(X_{n})\text{.}  \label{ind2}
\end{equation}%
Then $X=\Delta _{S}$ is defined by Formulas (\ref{ind2}) and (\ref{inters}).

For $X=C_{S}$ let $X_{0}=C_{0}$ be the square of diameter 1 in the Euclidean
quarter-plane $x\geq 0,y\geq 0$ with one side on each of the coordinate
axes. Denote by $v_{0},v_{1},v_{2},v_{3}$ its vertices and $%
w_{0},w_{1},w_{2},w_{3}$ the midpoints of its sides. Now define $X_{n}$
iteratively by the formulas 
\begin{equation}
X_{0}=C_{0},\quad X_{n+1}:=\bigcup_{i=0}^{3}(\psi _{\{\frac{1}{3}%
,v_{i}\}}(X_{n})\cup \psi _{\{\frac{1}{3},w_{i}\}}(X_{n}))\text{.}
\label{ind3}
\end{equation}%
Then $X=C_{S}$ is defined by Formulas (\ref{ind3}) and (\ref{inters}).

Denote by $\varepsilon _{n}$ the number $\frac{1}{2^{n+1}}$ in the cases $X=H
$ and $X=\Delta _{S}$ and the number $\frac{1}{2}\left( \frac{1}{3^{n}}%
\right) $ in the case $X=C_{S}$. Moreover, by our choice of $\varepsilon _{n}
$, $\phi _{X_{n}\varepsilon _{n}}:\left( X_{n}\right) _{\varepsilon
_{n}}\rightarrow X_{n}$ (see Notation \ref{metnote}) is the universal cover
of $X_{n}$, which means that the projection $\tau _{n}:\delta
_{1}(X_{n})\rightarrow \delta _{n}(X_{n}):=\delta _{\varepsilon _{n}}(X_{n})$
is an isomorphism. One can easily see that in all cases $X_{n}$ is a compact
Poincar\'{e} space and $\delta _{\varepsilon _{n}}(X)$ is isomorphic to $%
\delta _{\varepsilon _{n}}(X_{n})$ via the map $j_{n}=(i_{n})_{E_{%
\varepsilon _{n}}E_{\varepsilon _{n}}}$ defined in Definition \ref{induced},
where $i_{n}:X\rightarrow X_{n}$ is the inclusion. According to Corollary %
\ref{0518equiv}, $\lambda :\pi _{1}(X_{n})\rightarrow \delta _{1}(X_{n})$ is
an isomorphism. Therefore for any fixed $n$, the function $\omega _{n}:=\tau
_{n}\circ \lambda :\pi _{1}(X_{n})\rightarrow \delta _{n}(X_{n})$ is an
isomorphism.

According to Proposition \ref{lifter}, the image with respect to $\omega
_{n} $ of an equivalence class of a loop $c$ based at $\ast $ is the $E$%
-homotopy class of any sufficiently fine $E$-chain of the form $%
\{c(t_{0}),...,c(t_{k})\}$. Therefore, if $n\geq m$ we have the following
commutative diagram:%
\begin{equation*}
\begin{array}{lllll}
\pi _{1}(X_{n}) & \overset{\omega _{n}}{\longrightarrow } & \delta
_{n}(X_{n}) & \underrightarrow{^{j_{n}^{-1}}} & \delta _{\varepsilon _{n}}(X)
\\ 
\downarrow ^{\left( i_{mn}\right) _{\ast }} &  &  &  & \downarrow ^{\theta
_{mn}} \\ 
\pi _{1}(X_{m}) & \overset{\omega _{m}}{\longrightarrow } & \delta
_{m}(X_{m}) & \underrightarrow{^{j_{m}^{-1}}} & \delta _{\varepsilon _{m}}(X)%
\end{array}%
\end{equation*}%
where $i_{mn}:X_{n}\rightarrow X_{m}$ is inclusion and $j=h\circ \theta
_{\varepsilon _{m}\varepsilon _{n}}$ with $h:(X_{n})_{\varepsilon
_{m}}\rightarrow (X_{m})_{\varepsilon _{m}}$ being the inclusion-induced
mapping given by Definition \ref{induced}.

But $X_{n}$ is homotopic to a wedge product of $q_{n}$ circles, where $%
q_{n}=n$ for $X=H$, $q_{n}=\sum_{k=0}^{n-1}3^{k}$ for $X=\Delta _{S}$, and $%
q_{n}=\sum_{k=0}^{n-1}8^{k}$ for $X=C_{S}.$ So, in any case, $\pi
_{1}(X_{n})\cong F_{q_{n}},$ where $F_{q}$ is the free group with $q$
generators. Note that $\delta _{1}(X)=\underleftarrow{\lim }\delta
_{\varepsilon _{n}}(X)$ and for all three examples the inverse sequence $%
\left( \delta _{\varepsilon _{n}}(X),\theta _{mn}\right) $ is a cofinal
sequence in the inverse system $(F_{n},\pi _{mn})$ where $\pi
_{mn}:F_{n}\rightarrow F_{m}$ is the unique surjective homomorphism that
kills one extra generator.

Although all three spaces have the same deck group, the spaces $H,\Delta
_{S},$ and $C_{S}$ are mutually non-homeomorphic one-dimensional Peano
continua. In fact, the last two spaces are not semi-locally simply connected
at any point, while $H$ has unique point ($O$), at which it is not
semi-locally simply connected. The spaces $\Delta _{S}$ and $C_{S}$ are not
homeomorphic, because the first space have a countable subset $M$ such that $%
\Delta _{S}-M$ is null-homotopic in itself, while the second one has no such
subset. We can take 
\begin{equation*}
M=\bigcup_{n=0}^{\infty }\psi ^{n}(\{m\}),
\end{equation*}%
where $m$ is the midpoint of the segment $[v_{1},v_{2}],$ and 
\begin{equation*}
\psi (A):=\bigcup_{i=0}^{2}\psi _{\{\frac{1}{2},v_{i}\}}(A)
\end{equation*}%
for every subset $A\subset E^{n}.$

Since $\Delta _{S}$ and $C_{S}$ are non-homeomorphic one-dimensional Peano
continua that are not semi-locally simply connected at any point, then by a
result of \cite{Eda}, $\pi _{1}(\Delta _{S})$ and $\pi _{1}(C_{S})$ are not
isomorphic as abstract groups. Thus it follows from this and the above
result $\delta _{1}(\Delta _{S})\cong \delta _{1}(C_{S})$ that the
homomorphism $\pi _{1}(X)\rightarrow \delta _{1}(X)$ cannot be a bijection
for all three examples. Finally, note that the conditions of Proposition \ref%
{dense} are satisfied and hence the image of the fundamental group with
respect to $\lambda $ is dense in $\delta _{1}(X)$ in each case.

\begin{remark}
A similar discussion may be applied to the Sierpin'ski sponge $S_{S}$ (also
known as the universal Sierpin'ski curve) which is formed by successively
removing \textquotedblleft middle cores\textquotedblright\ from the unit
cube, implying that $\delta _{1}(S_{S})$ is also the inverse limit of
finitely generated free groups. It is well-known that every $1$-dimensional
locally connected metrizable continuum can be topologically embedded in $%
S_{S}$. By a result of Curtis and Fort (\cite{CF}), these embeddings induce
inclusions of the fundamental group, and it follows that the fundamental
group of any such continuum embeds in an inverse limit of finitely generated
free groups. This fact has been previously established by various authors,
most recently \cite{CC3}.
\end{remark}

\section{Topological Groups\label{topsec}}

The construction, in \cite{BPTG}, of the group $\widetilde{G}$ for a
(Hausdorff) topological group $G$ is the same as the construction that we
use in the present paper. However, the construction in the prior paper
includes a compatible group structure on each group in the fundamental
inverse system, which induces a group structure on $\widetilde{G}$ so that
the natural homomorphism $\phi :\widetilde{G}\rightarrow G$ is a quotient
map with closed, central, prodiscrete kernel. (There are also several
results in \cite{BPTG} that have no analogs for uniform spaces in general,
such as results concerning extensions of local homomorphisms.) For a
topological group $G$, the condition that we call \textquotedblleft
universal\textquotedblright\ in this paper is equivalent to what was called
\textquotedblleft locally defined\textquotedblright\ in \cite{BPTG}, by
Proposition 61 of \cite{BPTG}. A \textquotedblleft coverable
group\textquotedblright\ was defined to be the quotient of a locally defined
group via a closed normal subgroup. According to Corollary \ref{bi} in the
present paper, such a group is a coverable uniform space.

The converse of this statement involves Theorem 90 in \cite{BPTG}, which
requrires a correction. In fact, Professor Helge Gl\"{o}ckner of T. U.
Darmstadt has pointed out that the proof of Lemma 42 in \cite{BPTG} has a
gap, and we do not know whether this lemma is true. The only direct
reference to Lemma 42 is Theorem 90, and in light of this gap part (2) of
Theorem 90 should be restated as the following stronger condition: (2$%
^{\prime }$) $G$ has a basis for its topology at $e$ consisting of locally
generated symmetric neighborhoods, and $\phi _{U}$ is surjective for all $U$
in this basis. The proof may be modified to show (1) implies (2$^{\prime }$)
as follows. In the first part of the proof it is shown that $\phi $ is
surjective and that $G_{U}$ is coverable. The same argument may then be
applied to show that $\phi _{U}$ is surjective, proving (2$^{\prime }$). In
the proof of (2)$\Rightarrow $(3), Lemma 42 is only used to prove that all
of the homomorphisms $\phi _{U}$ are surjective, and so (2$^{\prime }$)
eliminates the need for Lemma 42. Note that (2$^{\prime }$) is actually a
stronger condition than the definition of coverable in the present paper,
which now completes the proof that a topological group is coverable in the
sense of \cite{BPTG} if and only if it is coverable in the present sense.

We know of no example of a uniform space $X$ (let alone topological group)
that is not coverable and $\phi :\widetilde{X}\rightarrow X$ is surjective,
hence it is still possible that Theorem 90 (or its generalization to uniform
spaces) is true as stated. However, the revised version is sufficient for
all applications in \cite{BPTG} and \cite{BPLCG} except for three. The first
is Corollary 91 of \cite{BPTG}, which is not used elsewhere. The second
exception is that the alternate hypothesis in Theorem 15 \textquotedblleft
or $\phi :\widetilde{G}\rightarrow G$ is surjective\textquotedblright\ needs
to be taken out (or replaced by a stronger assumption related to the new
condition (2$^{\prime }$)). The third exception is Corollary 107, which
isn't used elsewhere.

In the interest of completeness we will now address every reference to
Theorem 90 both in \cite{BPTG} and \cite{BPLCG} to show that the new version
is sufficient. Since (2$^{\prime }$) is stronger than (2), there is no
problem with any statement that doesn't involve (2)$\Rightarrow $(1) or (2)$%
\Rightarrow $(3). Theorem 90 is used this way in proofs of the following
statements in \cite{BPTG}: Theorem 4, Theorem 5, Proposition 10, Example 99,
Theorem 92 (necessity), Theorem 101, and Theorem 15 (second reference).

Now consider the remaining references: In Theorem 92 (sufficiency) of \cite%
{BPTG}, all of the projections are surjective and therefore (2$^{\prime }$)
holds. In \cite{BPLCG}, Theorem 24 is simply a restatement of Theorem 90
(with an additional statement about metric spaces added) and therefore (2)
must be replaced by (2$^{\prime }$). Theorem 24 is only used in the proof of
Theorem 7, and the only problematic usage is that Theorem 24, (2)$%
\Rightarrow $(3), is used to prove Theorem 7, (6)$\Rightarrow $(1). However,
it is well-known that an pathwise connected, locally compact group is
locally pathwise connected. Therefore for any pathwise connected symmetric
open set $U$ containing the origin, $G_{U}$ is pathwise connected and hence $%
\phi _{U}:\widetilde{G}\rightarrow G_{U}$ is surjective. That is, the
conditions for (2$^{\prime }$) are satisfied.

The relationship between covers in the present sense and covers in the sense
of \cite{BPTG} is considered in \cite{P}. Note that our paper partially
answers Problem 152 in \cite{BPTG}, which asks whether the generalized
fundamental group in that paper, which is the same as the deck group in the
present paper, is a topological invariant. From the present paper we know
that the deck group is in fact an invariant of uniform structures and for
compact groups a topological invariant. There are several other theorems
from \cite{BPTG} that likely can be generalized to the more setting of
uniform spaces, and these will be considered in a future paper. Some of the
questions of that paper have analogs for uniform spaces, most notably:

\begin{problem}
If $X$ is a uniform space, is $\widetilde{X}$ always universal?
\end{problem}

\begin{problem}
If $X$ is a uniform space and $\phi :\widetilde{X}\rightarrow X$ is a
uniform homeomorphism, is $X$ coverable?
\end{problem}

\begin{problem}
If $X$ is a simply connected, pathwise connected coverable space, is $X$
universal?
\end{problem}

\begin{acknowledgement}
The research for this paper was done while the first author enjoyed the
hospitality of the Department of Mathematics at the University of Tennessee
and was supported by RFBR grant N 05-01-00057-a.
\end{acknowledgement}

\end{document}